# BREAK DETECTION IN THE COVARIANCE STRUCTURE OF MULTIVARIATE TIME SERIES MODELS[1]

By Alexander Aue, Siegfried Hörmann, Lajos Horváth
and Matthew Reimherr

*University of California, Davis, University of Utah, University of Utah
and University of Chicago*

In this paper, we introduce an asymptotic test procedure to assess the stability of volatilities and cross-volatilites of linear and nonlinear multivariate time series models. The test is very flexible as it can be applied, for example, to many of the multivariate GARCH models established in the literature, and also works well in the case of high dimensionality of the underlying data. Since it is nonparametric, the procedure avoids the difficulties associated with parametric model selection, model fitting and parameter estimation. We provide the theoretical foundation for the test and demonstrate its applicability via a simulation study and an analysis of financial data. Extensions to multiple changes and the case of infinite fourth moments are also discussed.

**1. Introduction.** Univariate time series models, both linear and nonlinear, are as of today rather well understood. The literature furnishes a broad variety of contributions on the probabilistic structure of these models providing, for example, criteria for the existence of stationary solutions and the finiteness of moments. Manifold estimation procedures are well-developed and tests assessing the structural stability and the goodness of fit are widely available (see, among many others, [9, 16, 18, 39, 40, 48] and the references therein).

This, however, is no longer the case for multivariate time series models, in particular for multivariate financial time series. While many generalizations to higher-dimensional settings have been proposed since the inception of the

Received November 2008; revised March 2009.
[1]Supported in part by NSF Grants DMS-06-04670 and DMS-06-52420 and Grant RGC–HKUST6428/06H.

univariate ARCH and GARCH models by Engle [23] and Bollerslev [11] (see, for example, the survey articles [7, 47] and Section 4 below), probabilistic and statistical tools are still in developing stages. The underlying nontrivial theory is often times only rudimentary developed. Establishing efficient estimation procedures proves to be a challenging problem, since it bears the difficulty of incorporating a potentially large number of parameters into the algorithms. Tests that assess the structural stability of volatilities or cross-volatilities for multivariate nonlinear time series have yet to be introduced. On the other hand, multivariate nonlinear time series modeling is of high importance in practice as it is essential for an understanding of the coherence of exchange rates or market indexes. How movements of the financial markets are interpreted will significantly determine, for example, the asset allocation in a portfolio or the decision making of a risk manager (see [7, 24, 46, 47]). It is therefore crucial to understand the dynamic dependence structure in multivariate financial time series.

The main aim of this paper is to add to the research in this area by studying in detail the volatility and co-volatility structure of $d$-dimensional random vectors, which allow both linear and nonlinear specifications. We introduce an asymptotic test that is very flexible in that it only requires general but easily verifiable dependence assumptions to be satisfied by the underlying sequence. The procedure bypasses the difficulty of estimating a large number of model parameters and is very applicable also in high dimensions. As a byproduct, we provide stationarity conditions for several multivariate nonlinear time series models not previously established in the literature.

The results are based on an approximation of the underlying random vectors with a sequence of $m$-dependent random vectors which in turn yields a multivariate functional central limit theorem. This approximation generalizes similar univariate results provided in [10, 32] and is the main theoretical tool in the statistical analysis. In practice, this approach is generally more convenient to apply than competing methods such as the various notions of mixing (see, e.g., [17, 43]) as it typically requires less restrictive and easier verifiable assumptions. Moreover, the specifications employed here can even be modified to detect structural breaks in the co-movements of the $d$ coordinate processes under the assumptions of heavy tails (less than four finite moments) and even an infinite variance–covariance structure. This is especially important in applications related to financial time series.

An outline of the paper can be given as follows. In Section 2, we discuss the model assumptions and state the main results. These results are then extended to multiple change scenarios and to sequences with infinite variance–covariance structure in Section 3. Section 4 deals with a number of examples that are included in the provided framework. Section 5 is concerned with practical aspects and investigates the finite sample behavior through a simulation study and an application to stock index data. The multivariate functional central limit theorem is presented in Appendix A, while the mathematical proofs are given in Appendix B.

**2. Main results.** Let $(\mathbf{y}_j : j \in \mathbb{Z})$ be a sequence of $d$-dimensional random vectors with $E[\mathbf{y}_j] = \boldsymbol{\mu}$ and $E[|\mathbf{y}_j|^2] < \infty$. Here $|\cdot|$ denotes the Euclidean norm in $\mathbb{R}^d$. The main aim of the present paper is to introduce a testing procedure that identifies structural breaks in the volatilities and cross-volatilities of the process $(\mathbf{y}_j : j \in \mathbb{Z})$ based on observations of the random vectors $\mathbf{y}_1, \ldots, \mathbf{y}_n$. To do so, we consider the null hypothesis

$$(2.1) \qquad H_0 : \mathrm{Cov}(\mathbf{y}_1) = \cdots = \mathrm{Cov}(\mathbf{y}_n),$$

which indicates the constancy of the covariances in the observation period. It is usually under $H_0$ that a statistician can produce meaningful estimates and reliable forecasts. Whether or not $H_0$ holds should therefore precede any further statistical analysis. It should, however, be noted that a common approach to forecasting in the finance sector and elsewhere is to model volatilities via exponential smoothing with little regard to stationarity assumptions on the underlying processes. This works well as long as the volatilities can be assumed to evolve in a roughly smooth fashion, but leads to inaccurate forecasts in the nonsmooth break setting we consider in this paper.

As the alternative hypothesis we specify in particular the common scenario which allows for at most one change in the covariances. More precisely, it is assumed that there exists a—typically unknown—time lag $k^*$, referred to as the change-point, such that

$$(2.2) \quad H_A : \mathrm{Cov}(\mathbf{y}_1) = \cdots = \mathrm{Cov}(\mathbf{y}_{k^*}) \neq \mathrm{Cov}(\mathbf{y}_{k^*+1}) = \cdots = \mathrm{Cov}(\mathbf{y}_n).$$

Generalizations to more general alternatives allowing for several change-points are discussed briefly in Section 3 below. Since we are interested only in the second-order behavior of the random vectors $(\mathbf{y}_j : j \in \mathbb{Z})$, it is throughout this paper assumed that the expected values $E[\mathbf{y}_j]$ do not change over time. This may require an initial test of constancy of the means and, if necessary, a transformation of the data so that the expected values can be regarded as stable for the whole observation period. For further information on testing for changes in the mean, we refer to [18].

Even in the univariate case, there are relatively few contributions addressing changes in the variance. Gombay, Horváth and Hušková [26] and Inclán and Tiao [33] discuss the detection of (multiple) changes in the variance of independent observations based on weighted test statistics and an iterated cumulative sum of squares algorithm, respectively. The special case of (independent) Gamma distributions has been treated in [31]. More recently, [37] have introduced tests for a single change in the volatility of an ARCH($\infty$) series, while [20] analyze parametric nonlinear time series models by means of minimum description length procedures. There are, however, no in-depth treatments of (nonparametric) break detection methods for changes in the volatilities and cross-volatilities for multivariate, potentially nonlinear time series.

To construct a test statistic for distinguishing between $H_0$ and $H_A$, we let vech$(\cdot)$ be the operator that stacks the columns below the diagonal of a symmetric $d \times d$ matrix as a vector with $\mathfrak{d} = d(d+1)/2$ components. Assume for the moment that $E[\mathbf{y}_j] = \mathbf{0}$. Then, it holds under $H_0$ that the expected values $E[\text{vech}(\mathbf{y}_j \mathbf{y}_j^T)]$ are the same for $j = 1, \ldots, n$. Consequently, a version of the traditional cumulative sum (CUSUM) statistic can be constructed using the quantities

$$\mathcal{S}_k = \frac{1}{\sqrt{n}} \left( \sum_{j=1}^{k} \text{vech}[\mathbf{y}_j \mathbf{y}_j^T] - \frac{k}{n} \sum_{j=1}^{n} \text{vech}[\mathbf{y}_j \mathbf{y}_j^T] \right), \qquad k = 1, \ldots, n,$$

which basically compare the estimators of $E[\text{vech}(\mathbf{y}_j \mathbf{y}_j^T)]$ based on $k$ observation with the one computed from all observations. For $j = 1, \ldots, n$, let $\tilde{\mathbf{y}}_j = \mathbf{y}_j - \bar{\mathbf{y}}_n$ with $\bar{\mathbf{y}}_n = \frac{1}{n} \sum_{j=1}^{n} \mathbf{y}_j$. If $E[\mathbf{y}_j] \neq \mathbf{0}$, then the $\mathcal{S}_k$ can be replaced with the mean corrected modifications

$$\tilde{\mathcal{S}}_k = \frac{1}{\sqrt{n}} \left( \sum_{j=1}^{k} \text{vech}[\tilde{\mathbf{y}}_j \tilde{\mathbf{y}}_j^T] - \frac{k}{n} \sum_{j=1}^{n} \text{vech}[\tilde{\mathbf{y}}_j \tilde{\mathbf{y}}_j^T] \right), \qquad k = 1, \ldots, n.$$

Before we introduce the precise form of the test statistic, we detail the assumptions needed on the sequence $(\mathbf{y}_j : j \in \mathbb{Z})$. We are interested in several competing parametric specifications which are described in detail in Section 4 below. It is, however, often times hard to distinguish between different parametric models in practice and therefore assumptions imposed on $(\mathbf{y}_j : j \in \mathbb{Z})$ need to be general. Let $\alpha \geq 1$. For a random vector $\mathbf{x} \in \mathbb{R}^d$, let $|\mathbf{x}|$ be its Euclidean norm in $\mathbb{R}^d$ and $\|\mathbf{x}\|_\alpha = (E[|\mathbf{x}|^\alpha])^{1/\alpha}$ be its $L^\alpha$-norm.

ASSUMPTION 2.1. Let $(\mathbf{y}_j : j \in \mathbb{Z})$ be such that, under $H_0$, it satisfies the relations

$$(2.3) \qquad \mathbf{y}_j = \mathbf{f}(\varepsilon_j, \varepsilon_{j-1}, \ldots), \qquad j \in \mathbb{Z},$$

where $\mathbf{f} : \mathbb{R}^{d' \times \infty} \to \mathbb{R}^d$ is a measurable function and $(\varepsilon_j : j \in \mathbb{Z})$ a sequence of independent, identically distributed random variables with values in $\mathbb{R}^{d'}$. It is further required that, under $H_0$, there is a sequence of $m$-dependent random vectors $(\mathbf{y}_j^{(m)} : j \in \mathbb{Z})$ such that

$$(2.4) \qquad \mathbf{y}_j^{(m)} = \mathbf{f}^{(m)}(\varepsilon_j, \ldots, \varepsilon_{j-m}), \qquad j \in \mathbb{Z},$$

with measurable functions $\mathbf{f}^{(m)} : \mathbb{R}^{d' \times (m+1)} \to \mathbb{R}^d$, and

$$(2.5) \qquad \sum_{m \geq 1} \|\mathbf{y}_0 - \mathbf{y}_0^{(m)}\|_4 < \infty.$$

Condition (2.3) states in other words that $(\mathbf{y}_j : j \in \mathbb{Z})$ admits a causal representation, possibly nonlinear, in terms of the $d'$-dimensional innovation sequence $(\boldsymbol{\varepsilon}_j : j \in \mathbb{Z})$. A weak dependence structure is enabled through the introduction of the $m$-dependent random vectors $(\mathbf{y}_j^{(m)} : j \in \mathbb{Z})$ in (2.4) which are close to the original $(\mathbf{y}_j : j \in \mathbb{Z})$ in the sense of the distance measure (2.5). Note, however, that $(\mathbf{y}_j : j \in \mathbb{Z})$ itself need not be $m$-dependent. This type of condition has first been used in the context of univariate random variables in [32] and later also in [10, 29]. We will show in Section 4 that several popular multivariate time series models satisfy Assumption 2.1.

As shown in Appendix A, Assumption 2.1 induces a functional central limit theorem (FCLT) for the sequence $(\mathrm{vech}[\mathbf{y}_j \mathbf{y}_j^T] : j \in \mathbb{Z})$. FCLTs are frequently used in the asymptotic theory of weakly dependent processes. Various forms based, for example, on near epoch dependence and mixing conditions are discussed in [19, 43, 45]. The advantage of Assumption 2.1 over its mixing competitors is twofold. First, it is tailor-made for univariate and multivariate GARCH-type processes as the construction of $(\mathbf{y}_j^{(m)} : j \in \mathbb{Z})$ is virtually always straightforward (see [29] and Section 4 below), whereas mixing conditions are typically hard to verify and not applicable for several important time series models (see Section 4 below and [3, 29]). Second, mixing conditions require usually additional smoothness or restrictive moment assumptions (see [17]), while Assumption 2.1 gets by with finite fourth moments (which cannot be improved on when studying the covariance structure as in the present paper). These issues have led several authors such as [22, 50] to introduce more general dependence concepts which are more widely applicable. It is not immediately clear, however, whether the sharp almost sure invariance principles in [50] can be extended with reasonable effort to a multivariate setting. The methods in [22] are very general, but require more restrictive moment assumptions. In Section 3, we will address how to further generalize our methodology to multivariate sequences with infinite moments of less than fourth order.

Condition (2.5) yields in particular (see Theorem A.2 in Appendix A) that the long-run covariance

$$\Sigma = \sum_{j \in \mathbb{Z}} \mathrm{Cov}(\mathrm{vech}[\mathbf{y}_0 \mathbf{y}_0^T], \mathrm{vech}[\mathbf{y}_j \mathbf{y}_j^T])$$

converges (coordinatewise) absolutely. Statistical inference requires the estimation of $\Sigma$. We assume hence that there is an estimator $\hat{\Sigma}_n$ satisfying

$$(2.6) \qquad |\hat{\Sigma}_n - \Sigma|_E = o_P(1) \qquad (n \to \infty),$$

where, for a $d \times d$ matrix $M$, $|M|_E = \sup_{x \neq 0} |Mx|/|x|$ denotes the matrix norm induced by the Euclidean norm on $\mathbb{R}^d$. We shall discuss specific estimators $\hat{\Sigma}_n$ obeying the rate in (2.6) in Section 5. We now suggest the following two sets of test statistics to test $H_0$ against $H_A$. Define

$$\Lambda_n = \max_{1 \leq k \leq n} \mathcal{S}_k^T \hat{\Sigma}_n^{-1} \mathcal{S}_k \quad \text{and} \quad \Omega_n = \frac{1}{n} \sum_{k=1}^n \mathcal{S}_k^T \hat{\Sigma}_n^{-1} \mathcal{S}_k,$$

as well as

$$\tilde{\Lambda}_n = \max_{1 \leq k \leq n} \tilde{\mathcal{S}}_k^T \hat{\Sigma}_n^{-1} \tilde{\mathcal{S}}_k \quad \text{and} \quad \tilde{\Omega}_n = \frac{1}{n} \sum_{k=1}^n \tilde{\mathcal{S}}_k^T \hat{\Sigma}_n^{-1} \tilde{\mathcal{S}}_k.$$

The limit distributions of all four test statistics are given in the next theorem.

THEOREM 2.1 (Asymptotic under $H_0$). *Suppose that $(\mathbf{y}_j : j \in \mathbb{Z})$ are $d$-dimensional random vectors satisfying Assumption 2.1 and let $\hat{\Sigma}_n$ be such that (2.6) holds. Given that $E[\mathbf{y}_j] = \mathbf{0}$, it holds under (2.1) that*

$$(2.7) \qquad \Lambda_n \xrightarrow{\mathcal{D}} \Lambda(\mathfrak{d}) = \sup_{0 \leq t \leq 1} \sum_{\ell=1}^{\mathfrak{d}} B_\ell^2(t) \qquad (n \to \infty)$$

*and*

$$(2.8) \qquad \Omega_n \xrightarrow{\mathcal{D}} \Omega(\mathfrak{d}) = \sum_{\ell=1}^{\mathfrak{d}} \int_0^1 B_\ell^2(t) \, dt \qquad (n \to \infty),$$

*where $\mathfrak{d} = d(d+1)/2$, $(B_\ell(t) : t \in [0, 1])$, $1 \leq \ell \leq \mathfrak{d}$, are independent standard Brownian bridges and $\xrightarrow{\mathcal{D}}$ indicates convergence in distribution. The limit results in (2.7) and (2.8) remain valid if $\tilde{\Lambda}_n$ and $\tilde{\Omega}_n$ are used in place of $\Lambda_n$ and $\Omega_n$, respectively.*

The proof of Theorem 2.1 is relegated to Appendix B.1. The distribution of the random variables $\Lambda(\mathfrak{d})$ and $\Omega(\mathfrak{d})$ were derived in [36]. It is shown there that, for $x > 0$ and $\mathfrak{d} \geq 2$, the distribution function of $\Omega(\mathfrak{d})$ has the series expansion

$$(2.9) \begin{aligned} &P(\Omega(\mathfrak{d}) \leq x) \\ &= \frac{2^{(\mathfrak{d}+1)/2}}{\pi^{1/2} x^{\mathfrak{d}/4}} \sum_{j=0}^\infty \frac{\Gamma(j + \mathfrak{d}/2)}{j! \Gamma(\mathfrak{d}/2)} e^{-(j+\mathfrak{d}/4)^2/x} C_{(\mathfrak{d}-2)/2}\left(\frac{2j + \mathfrak{d}/2}{x^{1/2}}\right), \end{aligned}$$

where $\Gamma$ denotes the Gamma function and $C_{(\mathfrak{d}-2)/2}$ the parabolic cylinder functions (see, e.g., page 246 in [25]). Somewhat more involved formulas for the distribution of $\Lambda(\mathfrak{d})$ are also provided in [36]. Since $\mathfrak{d}$ could potentially be large, it is interesting to consider the distributions of $\Lambda(\mathfrak{d})$ and $\Omega(\mathfrak{d})$ as $\mathfrak{d} \to \infty$. We therefore define the standardized variables

$$\bar{\Lambda}(\mathfrak{d}) = \frac{\Lambda(\mathfrak{d}) - \mathfrak{d}/4}{(\mathfrak{d}/8)^{1/2}} \quad \text{and} \quad \bar{\Omega}(\mathfrak{d}) = \frac{\Omega(\mathfrak{d}) - \mathfrak{d}/6}{(\mathfrak{d}/45)^{1/2}}.$$

REMARK 2.1.   *If $\mathfrak{d} \to \infty$, then both $\bar{\Lambda}(\mathfrak{d}) \xrightarrow{\mathcal{D}} N(0,1)$ and $\bar{\Omega}(\mathfrak{d}) \xrightarrow{\mathcal{D}} N(0,1)$, where $N(0,1)$ denotes a standard normal random variable.*

The proof of Remark 2.1 may be found in Appendix B.1. Next, we turn our attention to the behavior of the test statistics under the alternative hypothesis $H_A$. To establish a consistency result, we need another set of assumptions which will be specified now.

ASSUMPTION 2.2.   Let $k^* = \lfloor \theta n \rfloor$ with some $\theta \in (0,1)$. Let $(\mathbf{y}_j : j \in \mathbb{Z})$ be such that the pre-change sequence $(\mathbf{y}_j : j \leq k^*)$ is strictly stationary and ergodic with $E[\mathbf{y}_0] = \mathbf{0}$ and $E[|\mathbf{y}_0|^2] < \infty$, and the post-change sequence $(\mathbf{y}_j : j > k^*)$ satisfies the relations

$$(2.10) \qquad\qquad \mathbf{y}_j = \mathbf{y}_j^* + \mathbf{z}_{j,n}, \qquad j \geq k^* + 1,$$

where $(\mathbf{y}_j^* : j \in \mathbb{Z})$ denotes a strictly stationary, ergodic sequence with $E[\mathbf{y}_0^*] = \mathbf{0}$ and $E[|\mathbf{y}_0^*|^2] < \infty$, and $(\mathbf{z}_{j,n} : j \in \mathbb{Z})$ is such that

$$(2.11) \qquad\qquad \max_{k* < j \leq n} |\mathbf{z}_{j,n}| = o_P(1) \qquad (n \to \infty).$$

As common in the literature, the time of change $k^*$ depends on the sample size $n$ (see [18]). This implies that we are actually dealing with an array under the alternative. For the sake of simplicity, this is suppressed in the notation. Roughly speaking, however, we are observing a stationary sequence until $k^*$. After the change, a new sequence starts with $\mathbf{y}_{k^*+1}$ as initial value. By assumption (2.11), this sequence is asymptotically stationary. The same condition (2.11) also ensures that the law of large numbers applies to the post-change sequence.

If the corresponding model parameter values change at time lag $k^*$, the conditions imposed by Assumption 2.2 are, in particular, satisfied for all of the specific parametric time series models discussed in Section 4 below. The following theorem establishes the consistency of the test procedures.

THEOREM 2.2 (Asymptotic under $H_A$).   *Suppose that $(\mathbf{y}_j : j \in \mathbb{Z})$ are $d$-dimensional random vectors satisfying Assumption 2.2. Then it holds under (2.2) that*

$$\sup_{t \in [0,1]} \left| \frac{1}{\sqrt{n}} \mathcal{S}_{\lfloor nt \rfloor} - \mathcal{S}^*(t) \right| = o_P(1) \qquad (n \to \infty),$$

*where*

$$\mathcal{S}^*(t) = \begin{cases} t(1-\theta)E[\mathrm{vech}[\mathbf{y}_0\mathbf{y}_0^T] - \mathrm{vech}[\mathbf{y}_0^*(\mathbf{y}_0^*)^T]], & t \in [0, \theta], \\ \theta(1-t)E[\mathrm{vech}[\mathbf{y}_0\mathbf{y}_0^T] - \mathrm{vech}[\mathbf{y}_0^*(\mathbf{y}_0^*)^T]], & t \in (\theta, 1]. \end{cases}$$

*With suitable modifications of the limit process* $(\mathcal{S}^*(t): t \in [0, 1])$, *a corresponding statement holds true if the process* $(\mathcal{S}_{\lfloor nt \rfloor}: t \in [0, 1])$ *is replaced with* $(\tilde{\mathcal{S}}_{\lfloor nt \rfloor}: t \in [0, 1])$.

Observe that if the mean of the $\mathrm{vech}[\mathbf{y}_j\mathbf{y}_j^T]$ changes, then the process $(\tilde{\mathcal{S}}_{\lfloor nt \rfloor}/\sqrt{n}: t \in [0, 1])$ will have a nonzero limit with probability one. This shows that our test has asymptotic power one. The proof of Theorem 2.2 is relegated to Appendix B.1. Once $H_0$ is rejected, a practitioner needs to locate the break point $k^*$ or, alternatively, the break point fraction (relative to the sample size) $\theta$. For this, we suggest the estimator

$$(2.12) \qquad \hat{\theta}_n = \frac{1}{n} \arg\max_{1 \le k \le n} \mathcal{S}_k^T \hat{\Sigma}_n^{-1} \mathcal{S}_k.$$

From Theorem 2.2, one can deduce intuitively that for large enough sample sizes $n$, the quadratic form $\mathcal{S}_{\lfloor nt \rfloor}^T \hat{\Sigma}_n^{-1} \mathcal{S}_{\lfloor nt \rfloor}$ will reach its maximum in the same argument as the limiting quadratic form $(\mathcal{S}^*(t))^T \Sigma^{-1} \mathcal{S}^*(t)$, thus yielding consistency. In fact, strong consistency of $\hat{\theta}_n$ can be proved under quite general conditions. Although a detailed discussion of this is beyond the scope of this paper, we will devote several paragraphs in Appendix B.2 to show that the relation

$$(2.13) \qquad |\hat{\theta}_n - \theta| = \mathcal{O}\left(\frac{\log\log n}{n}\right) \qquad \text{a.s.}$$

holds in many practical situations.

We continue with extensions of our methodology in Section 3 and examples in Section 4, while Section 5 examines practical aspects of the testing procedures via a simulation study and an application.

**3. Extensions.** We discuss in this section two directions of possible extensions of the theory introduced in Section 2. First, we show how our results may be applied also in the case of multiple changes in the underlying data. Second, we extend the approach to sequences of random vectors with heavy tails.

3.1. *Multiple break point detection.* In the literature of break detection schemes for univariate processes, a number of multiple change-point testing and estimation procedures are well established. They are, however, mainly concerned with mean changes (such as [6]), or require parametric model assumptions (such as [20]). The references for the univariate change in the variance case listed in the previous section have similar drawbacks. There are two papers by Andreou and Ghysels [1, 2] which address testing for multiple changes in the co-movements of financial assets from an empirical point of view without pursuing the corresponding asymptotic theory. This paper is thus novel in this direction. To incorporate the case of multiple breaks in the covariances into our setting, we can exchange $H_A$ with the following counterpart. Assume that, for some $r \geq 1$, there are time lags $1 < k_1^* < \cdots < k_r^* < n$ such that, for $s = 1, \ldots, r$,

$$H_A^* : \mathrm{Cov}(\mathbf{y}_{k_{s-1}^*+1}) = \cdots = \mathrm{Cov}(\mathbf{y}_{k_s^*}) \neq \mathrm{Cov}(\mathbf{y}_{k_s^*+1}) = \cdots$$
$$= \mathrm{Cov}(\mathbf{y}_{k_{s+1}^*}),$$

where $k_0^* = 0$ and $k_{r+1}^* = n$. If the change locations have the standard form used in the analysis of multivariate changes and if the change size is either constant or shrinking at an appropriate rate (see [6, 18]), then, all four test statistics $\Lambda_n$, $\Omega_n$, $\tilde{\Lambda}_n$ and $\tilde{\Omega}_n$ have asymptotic power one. To make this statement precise, we replace Assumption 2.2 with the following requirements.

ASSUMPTION 3.1. Let $k_s^* = \lfloor n\theta_s \rfloor$, $s = 0, \ldots, r+1$, with $0 = \theta_0 < \theta_1 < \cdots < \theta_r < \theta_{r+1} = 1$. Let the sequence $(\mathbf{y}_j : j = 1, \ldots, n)$ be such that $\mathbf{y}_j = \mathbf{y}_{j,0}^*$, $j = 1, \ldots, k_1^*$ and

$$\mathbf{y}_j = \mathbf{y}_{j,s}^* + z_{j,n}^{(s)}, \qquad j = k_s^*, \ldots, k_{s+1}^* - 1, s = 1, \ldots, r,$$

where, for $s = 0, \ldots, r$, $(y_{j,s}^* : j \in \mathbb{Z})$ denotes a strictly stationary, ergodic sequence with $E[\mathbf{y}_{0,s}^*] = 0$ and $E[|\mathbf{y}_{0,s}|^2] < \infty$, and, for $s = 1, \ldots, r$, $(z_{j,n}^{(s)} : j \in \mathbb{Z})$ is such that

$$\max_{1 \leq s \leq r} \max_{k_s^* < j \leq k_{s+1}^*} |z_{j,n}^{(s)}| = o_P(1) \qquad (n \to \infty).$$

With Assumption 3.1, we obtain the following counterpart of Theorem 2.2 under the more general alternative $H_A^*$. A simple computation shows that the one-change case is reproduced if $r = 1$ is chosen.

THEOREM 3.1 (Asymptotic under $H_A^*$). *Suppose that $(\mathbf{y}_j : j = 1, \ldots, n)$ are $d$-dimensional random vectors satisfying Assumption 3.1. Then it holds under $H_A^*$ that*

$$\sup_{t \in [0,1]} \left| \frac{1}{\sqrt{n}} \mathcal{S}_{\lfloor nt \rfloor} - S_r^*(t) \right| = o_P(1) \qquad (n \to \infty),$$

*where* $S_r^*(0) = 0$ *and*

$$\mathcal{S}_r^*(t) = \sum_{u=0}^{r-1} (\theta_{u+1} - \theta_u) E[\mathrm{vech}[\mathbf{y}_{0,u}^*(\mathbf{y}_{0,u}^*)^T]] + (t - \theta_s) E[\mathrm{vech}[\mathbf{y}_{0,s}^*(\mathbf{y}_{0,s}^*)^T]] - t\vartheta$$

*for* $t \in (\theta_s, \theta_{s+1}]$, $s = 0, \ldots, r$, *with*

$$\vartheta = \sum_{s=0}^{r} (\theta_{s+1} - \theta_s) E[\mathrm{vech}[\mathbf{y}_{0,s}^*(\mathbf{y}_{0,s}^*)^T]].$$

*With suitable modifications of the limit process* $(\mathcal{S}_r^*(t) : t \in [0, 1])$, *a corresponding statement holds true if the process* $(\mathcal{S}_{\lfloor nt \rfloor} : t \in [0, 1])$ *is replaced with* $(\bar{\mathcal{S}}_{\lfloor nt \rfloor} : t \in [0, 1])$.

One can now easily verify that $\sup_{t \in [0,1]} |\mathcal{S}_r^*(t)| > 0$ under $H_A^*$, so that all four test statistics have asymptotic power one. In the empirical Section 5, we illustrate that the binary segmentation procedure based on $\Omega_n$ produces very reasonable results in an application to financial data and is therefore capable to detect multiple breaks in the variance–covariance structure of multivariate financial time series. Moreover, applying the technique developed in [6] would also imply the weak consistency of the estimators of $\theta_1, \ldots, \theta_r$, which are defined analogously as in (2.12) and are evaluated utilizing the binary segmentation procedure.

### 3.2. *Transformations for sequences with less than four finite moments.*
There is plenty of evidence in the literature that the assumption of four finite moments for financial time series is rather restrictive in applications. Theoretical results such as limit theorems for the sample autocorrelation function and the quasi-likelihood estimation of the unknown model parameters has for heavy-tailed univariate ARCH and GARCH processes been discussed in [21, 27, 41], among others. Even in the univariate case, these results are generally difficult to apply in practice because the statistical inference procedures require, for example, precise knowledge of the observations' number of finite moments. We are unaware of any results concerned with the corresponding multivariate versions of heavy-tailed GARCH processes. The approach introduced in Section 2, on the other hand, is in principal not hampered by these limitations as the following argument shows. Let again $(\mathbf{y}_j : j \in \mathbb{Z})$ denote the multivariate process to be observed. For $\delta \in (0, 1]$, let $(\mathbf{x}_j : j \in \mathbb{Z})$ be the fractional process transformation given by

$$(3.1) \qquad \mathbf{x}_j = |\mathbf{y}_j|^\delta = (|y_{j,1}|^\delta, \ldots, |y_{j,d}|^\delta)^T, \qquad j \in \mathbb{Z}.$$

Instead of detecting breaks in the covariance structure of the original $(\mathbf{y}_j : j \in \mathbb{Z})$, the transformed process $(\mathbf{x}_j : j \in \mathbb{Z})$ can be utilized. This bears the further advantage of increased flexibility of the testing procedure as it is in the transformed version also possible to find, for example, breaks in the marginal distribution of the $\mathbf{y}_j$. To apply our theory to (3.1), condition (2.5) needs to be verified for the $\mathbf{x}_j$. Let therefore $\mathbf{x}_j^{(m)} = |\mathbf{y}_j^{(m)}|^\delta$ be the sequence of $m$-dependent variables used for the approximation. Then $\|\mathbf{x}_j - \mathbf{x}_j^{(m)}\|_4 \leq d^{1/4} \|\mathbf{y}_j - \mathbf{y}_j^{(m)}\|_{4\delta}^\delta$ and consequently

$$(3.2) \qquad \sum_{m \geq 1} \|\mathbf{x}_0 - \mathbf{x}_0^{(m)}\|_4 \leq d^{1/4} \sum_{m \geq 1} \|\mathbf{y}_0 - \mathbf{y}_0^{(m)}\|_{4\delta}^\delta < \infty$$

provides us with a sufficient criterion. Picking $\delta < 1/2$ shows moreover that our method works theoretically even if the process $(\mathbf{y}_j : j \in \mathbb{Z})$ does not possess a finite variance–covariance structure. In the case of GARCH sequences, for example, it appears convenient to pick $\delta = 1/2$. In this case, (3.2) is even easier to verify than the original version in (2.5). We further discuss the choice $\delta < 1/2$ for the constant conditional correlation GARCH model in Section 4.2 below. It would be worthwhile developing a data-driven procedure to determine an optimal value of $\delta$. Heavy-tailed random vectors naturally display a greater amount of variability compared to light-tailed ones, so that it is generally much harder to discriminate between wild but stationary fluctuations and an actual break in the underlying process structure. On the other hand, choosing $\delta < 1$ usually flattens the signal, so that breaks may be harder to detect. A selection procedure for $\delta$ needs to balance these two properties. We leave this issue to future research.

**4. Examples.** In this section, we discuss the leading examples included in the testing framework developed in Section 2. As a special case, our theory applies to the most popular multivariate linear time series models such as vector-valued ARMA and linear processes (see Section 4.1). Our focus, however, is more on multivariate nonlinear financial time series. As examples, we introduce here Bollerslev's [12] constant conditional correlation GARCH model (Section 4.2), an extension of this model introduced by Jeantheau [34] (Section 4.3), factor GARCH models originally considered in Engle, Ng and Rothschild [24] (Section 4.4) and multivariate exponential GARCH models (Section 4.5).

4.1. *Multivariate ARMA and linear processes.* The theory of linear multivariate time series is closely connected to investigating the properties of linear processes. Let $(\boldsymbol{\varepsilon}_j : j \in \mathbb{Z})$ be a sequence of $d$-dimensional independent and identically distributed random variables with mean $\mathbf{0}$ and covariance matrix $\Psi$. A $d$-dimensional linear process $(\mathbf{y}_j : j \in \mathbb{Z})$ can then be constructed by letting

$$(4.1) \qquad \mathbf{y}_j = \sum_{\ell=0}^\infty C_\ell \boldsymbol{\varepsilon}_{j-\ell}, \qquad j \in \mathbb{Z},$$

where $(C_\ell : \ell \geq 0)$ is a sequence of $d \times d$ matrices with absolutely summable components. It is clear that $(\mathbf{y}_j : j \in \mathbb{Z})$ is a strictly stationary but dependent process that can be included in the setting of Section 2.

THEOREM 4.1. *A multivariate linear process* $(\mathbf{y}_j : j \in \mathbb{Z})$ *specified by* (4.1) *satisfies Assumption 2.1 if* $E[|\boldsymbol{\varepsilon}_0|^4] < \infty$ *and*

$$\sum_{m \geq 1} \sum_{\ell \geq m+1} |C_\ell| < \infty.$$

The proof of Theorem 4.1 is given in Appendix B. Arguably, the most useful class of multivariate linear processes is given by ARMA$(p, q)$ models. These are defined in terms of the linear difference equations

$$(4.2) \qquad \mathbf{y}_j - \sum_{\ell=1}^{p} A_\ell \mathbf{y}_{j-\ell} = \boldsymbol{\varepsilon}_j + \sum_{\ell=1}^{q} B_\ell \boldsymbol{\varepsilon}_{j-\ell}, \qquad j \in \mathbb{Z},$$

where $A_1, \ldots, A_p$ and $B_1, \ldots, B_q$ are $d \times d$ matrices. The following corollary identifies those multivariate ARMA models that are special cases of the linear process in (4.1).

COROLLARY 4.1. *A multivariate* ARMA$(p, q)$ *process* $(\mathbf{y}_j : j \in \mathbb{Z})$ *specified by* (4.2) *satisfies Assumption 2.1 if* $E[|\boldsymbol{\varepsilon}_0|^4] < \infty$ *and* $\det A(z) \neq 0$ *for all* $z \in \mathbb{C}$ *such that* $|z| \leq 1$, *where* $A(z) = 1 - A_1 z - \cdots - A_p z^p$ *is a matrix-valued polynomial.*

PROOF. It follows from Theorem 11.3.1 in [16] that the condition $\det A(z) \neq 0$ implies the existence of a unique strictly stationary solution to the difference equations (4.2) which is given in the form (4.1). Moreover, the $(C_\ell : \ell \geq 0)$ in (4.1) are uniquely determined by the matrix-valued power series

$$C(z) = \sum_{\ell=0}^{\infty} C_\ell z^\ell = A^{-1}(z) B(z), \qquad |z| \leq 1,$$

and decay at a geometric rate, whence we have verified the summability condition in Theorem 4.1 and complete the proof. □

4.2. *The constant conditional correlation GARCH model.* In this subsection, we discuss one possible way of specifying multivariate models that exhibit heteroscedasticity. It is a nontrivial task to extend the univariate ARCH and GARCH settings as introduced by Engle [23] and Bollerslev [11] to the vector case. Bollerslev [12] suggested the following constant conditional correlation (CCC) model. Let $(\boldsymbol{\varepsilon}_j : j \in \mathbb{Z})$ be a sequence of $d$-dimensional independent and identically distributed random vectors with mean $\mathbf{0}$ and covariance matrix $\Psi$. Let the process $(\mathbf{y}_j : j \in \mathbb{Z})$ be defined by

$$(4.3) \qquad \mathbf{y}_j = \boldsymbol{\sigma}_j \circ \boldsymbol{\varepsilon}_j,$$

$$(4.4) \qquad \boldsymbol{\sigma}_j \circ \boldsymbol{\sigma}_j = \boldsymbol{\omega} + \sum_{\ell=1}^{p} \boldsymbol{\alpha}_\ell \circ \boldsymbol{\sigma}_{j-\ell} \circ \boldsymbol{\sigma}_{j-\ell} + \sum_{\ell=1}^{q} \boldsymbol{\beta}_\ell \circ \mathbf{y}_{j-\ell} \circ \mathbf{y}_{j-\ell},$$

where $\boldsymbol{\omega}$ is coordinatewise strictly positive, $\boldsymbol{\alpha}_1, \ldots, \boldsymbol{\alpha}_p$ and $\boldsymbol{\beta}_1, \ldots, \boldsymbol{\beta}_q$ are coordinatewise nonnegative $d$-dimensional vectors, and $p$ and $q$ are nonnegative integers. Moreover, $\circ$ denotes the Hadamard product of two identically sized vectors (or, matrices), which is computed by elementwise multiplication. A process given by (4.3) and (4.4) formally resembles the structure of a univariate GARCH$(p, q)$ time series. In fact, each coordinate is specified by a one-dimensional GARCH equation, whose orders are at most $(p, q)$. (Lower orders can easily be achieved through zero coefficients.)

THEOREM 4.2. *A multivariate CCC–GARCH process* $(\mathbf{y}_j : j \in \mathbb{Z})$ *specified by (4.3) and (4.4) satisfies Assumption 2.1 if* $E[|\varepsilon_0|^4] < \infty$ *and*

$$\gamma_C = \max_{1 \le i \le d} \sum_{\ell=1}^{r} \|\alpha_{\ell,i} + \beta_{\ell,i} \varepsilon_{0,i}^2\|_2 < 1,$$

*where* $r = \max\{p, q\}$, $\boldsymbol{\alpha}_\ell = (\alpha_{\ell,1}, \ldots, \alpha_{\ell,d})^T$ *and* $\boldsymbol{\beta}_\ell = (\beta_{\ell,1}, \ldots, \beta_{\ell,d})^T$. *Additional coefficients appearing in the latter display are set to equal zero.*

The proof of Theorem 4.2 is given in Appendix B. The main condition, $\gamma_C < 1$, ensures that the process $(\mathbf{y}_j : j \in \mathbb{Z})$ has finite fourth-order moments and is easy to verify based on a decomposition given in Appendix B. Nelson [42] showed that $\gamma_C < 1$ is also necessary in the case of scalar GARCH$(1,1)$ models. Note, however, that even though necessary and sufficient conditions for the existence of fourth-order moments are established in the literature (see He and Teräsvirta [28]), it is in general difficult to decide whether a given parameterization satisfies these requirements.

To further demonstrate the applicability of our method to the transformed process $(\mathbf{x}_j : j \in \mathbb{Z})$ defined in (3.1), we show now that (3.2) can be verified for the multivariate CCC–GARCH model even if $\delta < 1/2$. For simplicity of the presentation, we assume for the moment that $p = q = 1$ and that each coordinate process $y_{j,i}$ possesses $4\delta$ moments ($\delta < 1/2$). Then, Theorem 3 of [42] implies that

$$\tilde{\gamma}_C = \max_{1 \le i \le d} E[(\alpha_i + \beta_i \varepsilon_{0,i}^2)^{2\delta}] < 1$$

is necessary and sufficient for $E[|y_{0,i}|^{4\delta}] < \infty$. Moreover, using the GARCH$(1,1)$ representation obtained in [8] (see also [42] and display (B.4) in Appendix B) and the truncated variables $\mathbf{y}_j^{(m)}$ in (B.5), one finds that, for all $i = 1, \ldots, d$,

$$
\begin{aligned}
\|y_{0,i} - y_{0,i}^{(m)}\|_{4\delta}^{\delta} &= \left( |\omega_i| E[|\varepsilon_{0,i}|^{4\delta}] E\left[ \left( \sum_{n=m+1}^{\infty} \prod_{k=1}^{n} (\alpha_i + \beta_i \varepsilon_{-k,i}^2) \right)^{2\delta} \right] \right)^{1/4} \\
&\leq \left( |\omega_i| E[|\varepsilon_{0,i}|^{4\delta}] E\left[ \sum_{n=m+1}^{\infty} \prod_{k=1}^{n} (\alpha_i + \beta_i \varepsilon_{-k,i}^2)^{2\delta} \right] \right)^{1/4} \\
&= \left( |\omega_i| E[|\varepsilon_{0,i}|^{4\delta}] \sum_{n=m+1}^{\infty} (E[(\alpha_i + \beta_i \varepsilon_{0,i}^2)^{2\delta}])^n \right)^{1/4} \\
&= \mathcal{O}(\tilde{\gamma}_C^{m/4}),
\end{aligned}
$$

using the independence of the $(\varepsilon_{j,i} \colon j \in \mathbb{Z})$ and the relation $|\sqrt{a} - \sqrt{b}|^2 \leq |a - b|$ for all $a, b \geq 0$ to obtain the first equality sign, and $2\delta < 1$ for the first inequality sign. The preceding proves via (3.2) that (2.5) holds for the transformed variables $\mathbf{x}_j = |\mathbf{y}_j|^{\delta}$ if $\delta < 1/2$. To employ the break detection procedure of Section 2 consequently requires not even a finite variance–covariance structure if dealing with a multivariate CCC–GARCH process. A similar argument applies also in the finite second moment case $\delta \geq 1/2$.

We conclude this subsection with the remark that the traditional GARCH specification in (4.4) can easily be replaced with other specifications that take into account, for example, asymmetries, as long as necessary and sufficient conditions for the existence of strictly stationary solutions and moments are available. One such class of GARCH models has been studied in Aue, Berkes and Horváth [5].

4.3. *Jeantheau's constant conditional correlation GARCH model.* An extension of the CCC–GARCH model was introduced by Jeantheau [34] by replacing the coefficients $\boldsymbol{\alpha}_1, \ldots, \boldsymbol{\alpha}_p$ and $\boldsymbol{\beta}_1, \ldots, \boldsymbol{\beta}_q$ in (4.4) with entrywise nonnegative $d \times d$ matrices $A_1, \ldots, A_p$ and $B_1, \ldots, B_q$. This leads to the modified definition

$$
(4.5) \qquad \boldsymbol{\sigma}_j \circ \boldsymbol{\sigma}_j = \boldsymbol{\omega} + \sum_{\ell=1}^{p} A_\ell (\boldsymbol{\sigma}_{j-\ell} \circ \boldsymbol{\sigma}_{j-\ell}) + \sum_{\ell=1}^{q} B_\ell (\mathbf{y}_{j-\ell} \circ \mathbf{y}_{j-\ell})
$$

for the conditional covariances. Therein, $\boldsymbol{\omega}$ is coordinatewise strictly positive. Letting $A_\ell = \operatorname{diag}(\boldsymbol{\alpha}_\ell)$ and $B_\ell = \operatorname{diag}(\boldsymbol{\beta}_\ell)$, obviously returns (4.4). In this case, the criteria for the existence of a unique, nonanticipative (i.e., future independent) strictly stationary solution can be obtained mimicking the univariate results provided by Brandt [15] and Bougerol and Picard [13, 14].

The more general case (4.5) does not follow straightforwardly and requires a number of modifications.

We introduce more notation. For $\alpha \in (0,1)$, we extend the $L^\alpha$-norm by letting here $\|\cdot\|_\alpha = E[|\cdot|^\alpha]$. For $\alpha \geq 1$, we continue to work with $\|\cdot\|_\alpha = (E[|\cdot|^\alpha])^{1/\alpha}$. Notice that the matrix norm $|\cdot|_E$ defined below display (2.6) satisfies the condition $|MN|_E \leq |M|_E |N|_E$ for two $d \times d$ matrices $M$ and $N$. For a $d \times d$ matrix $M$, let finally

$$\|M\|_{E,\alpha} = \||M|_E\|_\alpha$$

be the composition of the matrix norm $|\cdot|_E$ and the $L^\alpha$-norm $\|\cdot\|_\alpha$. Now, we can give a sufficient condition for the existence of a strictly stationary solution in Jeantheau's model.

THEOREM 4.3. *If $E[|\boldsymbol{\varepsilon}_0|^{2\alpha}] < \infty$ and*

$$\gamma_{J,\alpha} = \sum_{\ell=1}^r \|A_\ell + B_\ell E_0\|_{E,\alpha} < 1$$

*for some $\alpha > 0$, then the (4.5) have the unique, nonanticipative strictly stationary and ergodic solution*

$$\boldsymbol{\sigma}_j \circ \boldsymbol{\sigma}_j = \boldsymbol{\omega} + \left[\sum_{k=1}^\infty \sum_{1 \leq i_1 < \cdots < i_k \leq r} \prod_{\ell=1}^k (A_{i_\ell} + B_{i_\ell} E_{j-i_1-\cdots-i_\ell})\right] \boldsymbol{\omega}, \qquad j \in \mathbb{Z},$$

*where $r = \max\{p,q\}$ and $E_j = \operatorname{diag}(\boldsymbol{\varepsilon}_j \circ \boldsymbol{\varepsilon}_j)$. Additional matrices appearing in the definition of $\gamma_{J,\alpha}$ are set to equal the zero matrix.*

With this result at hand, it can further be established under which conditions the testing procedures of Section 2 are applicable for Jeantheau's model. The proofs of Theorems 4.3 and 4.4 are given in Appendix B.

THEOREM 4.4. *A multivariate CCC–GARCH process $(\mathbf{y}_j : j \in \mathbb{Z})$ specified by (4.3) and (4.5) satisfies Assumption 2.1 if $E[|\boldsymbol{\varepsilon}_0|^4] < \infty$ and $\gamma_J = \gamma_{J,2} < 1$.*

Observe, however, that $\gamma_J < 1$ here and $\gamma_C < 1$ in Bollerslev's [12] CCC–GARCH model are not even equivalent if $A_\ell = \operatorname{diag}(\boldsymbol{\alpha}_\ell)$ and $B_\ell = \operatorname{diag}(\boldsymbol{\beta}_\ell)$ are diagonal matrices. Rather, the condition $\gamma_C < 1$ is less restrictive. This can be seen from the inequality

$$\gamma_J = \sum_{\ell=1}^r \|A_\ell + B_\ell E_0\|_{E,2} \geq \max_{1 \leq i \leq d} \sum_{\ell=1}^r \|\alpha_{\ell,i} + \beta_{\ell,i} \varepsilon_{0,i}^2\|_2 = \gamma_C.$$

Standard arguments imply that this inequality is strict if, for example, $\boldsymbol{\varepsilon}_0$ is normally distributed (or, if it has a positive density on $\mathbb{R}^d$).

4.4. *Dynamic factor models.* It is often believed in economic theory that the dynamics of a multivariate sequence $(\mathbf{y}_j : j \in \mathbb{Z})$ can be adequately described by so-called (unobserved) common factors, say, $(\mathbf{z}_j : j \in \mathbb{Z})$ with values in $\mathbb{R}^{d'}$, $d' < d$. This gives rise to the model

$$(4.6) \qquad \mathbf{y}_j = L\mathbf{z}_j + \boldsymbol{\xi}_j, \qquad j \in \mathbb{Z},$$

where $L$ is a $d \times d'$ matrix of factor loadings and $(\boldsymbol{\xi}_j : j \in \mathbb{Z})$ a sequence of $d'$-dimensional errors. The following theorem specifies conditions under which the sequence $(\mathbf{y}_j : j \in \mathbb{Z})$ satisfies Assumption 2.1.

THEOREM 4.5. *Let $(\mathbf{z}_j : j \in \mathbb{Z})$ and $(\boldsymbol{\xi}_j : j \in \mathbb{Z})$ be such that they satisfy Assumption 2.1. Then, the process $(\mathbf{y}_j : j \in \mathbb{Z})$ specified by (4.6) satisfies Assumption 2.1 as well.*

The main argument for the use of factor models is their potential to reduce dimensionality and consequently the number of parameters to be estimated from data. Often, the components of the factors are assumed to be uncorrelated, indicating genuinely different driving forces that account for the overall volatility behavior. The first factor GARCH model was introduced by Engle, Ng and Rothschild [24], more recent factor models were put forth in [38, 49]. For an in-depth overview we refer to [47].

4.5. *Multivariate exponential GARCH.* A further subclass of the multivariate GARCH family is provided by Kawakatsu's [35] matrix extension of Nelson's [42] exponential GARCH model. The model in [35] is given by first assuming that the random vectors $(\mathbf{y}_j : j \in \mathbb{Z})$ follow the equations

$$(4.7) \qquad \mathbf{y}_j = H_j^{1/2} \boldsymbol{\varepsilon}_j, \qquad j \in \mathbb{Z},$$

with $(\boldsymbol{\varepsilon}_j : j \in \mathbb{Z})$ being a sequence of independent, identically distributed $d$-dimensional random vectors. The $d \times d$ matrices $H_j$ are now specified by letting

$$
\begin{aligned}
(4.8) \quad \text{vech}[\log H_j - C] &= \sum_{\ell=1}^{p} A_\ell \, \text{vech}[\log H_{j-\ell} - C] \\
&\quad + \sum_{\ell=1}^{q} B_\ell \boldsymbol{\varepsilon}_{j-\ell} + \sum_{\ell=1}^{q} F_\ell (|\boldsymbol{\varepsilon}_{j-\ell}|_c - E[|\boldsymbol{\varepsilon}_{j-\ell}|_c]),
\end{aligned}
$$

where $C$ denotes a symmetric $d \times d$ matrix and $|\cdot|_c$ component-wise absolute value. The parameter matrices $A_\ell$, $B_\ell$ and $F_\ell$ are of dimensions $\mathfrak{d} \times \mathfrak{d}$, $\mathfrak{d} \times d$ and $\mathfrak{d} \times d$, respectively. Since $\log H_j$ is clearly a symmetric matrix, we have that

$$(4.9) \qquad H_j = \exp(\log H_j) = \sum_{k=0}^{\infty} \frac{1}{k!} (\log H_j)^k, \qquad j \in \mathbb{Z},$$

is a positive definite matrix. The apparent advantage of this approach is that no (unnatural) restrictions on the parameters have to be imposed through the model specification so as to ensure positive definiteness. The dynamic correlation structure is enabled through the recursions in (4.8) which determine $\log H_j$, while $H_j$ in turn is obtained from the matrix exponential representation in (4.9). We shall work here with a slight modification given by

$$(4.10) \qquad \text{vech}[\log H_j - C] = A \, \text{vech}[\log H_{j-1} - C] + F(\varepsilon_{j-1}, \ldots, \varepsilon_{j-q}),$$

where $A$ is a $\mathfrak{d} \times \mathfrak{d}$ matrix and $F$ is a measurable function taking values in $\mathbb{R}^{\mathfrak{d}}$. The following theorem provides a sufficient condition for the existence of a strictly stationary solution to the multivariate exponential GARCH equations. Denote by $\log^+ x = \max\{\log x, 0\}$.

THEOREM 4.6. *If* $E[\log^+ |F(\varepsilon_q, \ldots, \varepsilon_1)|] < \infty$ *and*

$$(4.11) \qquad \limsup_{k \to \infty} \frac{1}{k} \log |A^k|_E < 0,$$

*then the (4.10) have the unique, nonanticipative strictly stationary and ergodic solution*

$$\text{vech}[\log H_j - C] = \sum_{k=0}^{\infty} A^k F(\varepsilon_{j-k-1}, \ldots, \varepsilon_{j-k-q}), \qquad j \in \mathbb{Z}.$$

*Condition (4.11) is in particular satisfied if* $|A|_E < 1$.

To prove that Assumption 2.1 is applicable here, we need to assume additionally that the random variable $|F(\varepsilon_{-1}, \ldots, \varepsilon_{-q})|$ has a finite moment generating function in a sufficiently large neighborhood of zero. Details are provided next.

THEOREM 4.7. *A multivariate exponential GARCH process* $(\mathbf{y}_j : j \in \mathbb{Z})$ *specified by (4.7) and (4.10) satisfies Assumption 2.1 if* $|A|_E < 1$ *and* $E[\exp(t|F(\varepsilon_{-1}, \ldots, \varepsilon_{-q})|)] < \infty$ *for some* $t > \sqrt{8}q$.

The proofs of Theorems 4.6 and 4.7 are postponed to Section B.4.

4.6. *Summary.* In this section, we have discussed the most popular multivariate GARCH models. Which of these models a practitioner will eventually pick will depend heavily on the application at hand. All specifications have their particular strengths and weaknesses. Research—applied and theoretic—is still to be conducted to reveal further insight. We refrain hence from giving recommendations here but refer to the survey papers [7, 47] for additional reading. Independent of practical modeling considerations, the tests introduced in Section 2, however, are flexible enough to cover any of the particular multivariate specifications given in Sections 4.1–4.5. To investigate their finite sample behavior will be subject of Section 5.

Table 1
*Critical values for $\bar{\Omega}(\mathfrak{d})$ for the three significance levels 0.10, 0.05 and 0.01. For each critical value, we report $P(\bar{\Omega}(\mathfrak{d}) \leq z_{1-\alpha})$ in brackets. The last column, $\mathfrak{d} = \infty$, corresponds to the standard normal case as stated in Remark 2.1*

| $\mathfrak{d}$ | 10 | 15 | 20 | 50 | 100 | 200 | 500 | $\infty$ |
|---|---|---|---|---|---|---|---|---|
| $q_{0.90}^{\Omega}(\mathfrak{d})$ | 1.33 (0.89) | 1.33 (0.89) | 1.32 (0.89) | 1.31 (0.89) | 1.31 (0.90) | 1.30 (0.90) | 1.29 (0.90) | 1.28 (0.90) |
| $q_{0.95}^{\Omega}(\mathfrak{d})$ | 1.84 (0.93) | 1.81 (0.94) | 1.79 (0.94) | 1.74 (0.94) | 1.71 (0.94) | 1.69 (0.94) | 1.68 (0.95) | 1.64 (0.95) |
| $q_{0.99}^{\Omega}(\mathfrak{d})$ | 2.90 (0.98) | 2.80 (0.98) | 2.74 (0.98) | 2.59 (0.98) | 2.51 (0.99) | 2.46 (0.99) | 2.41 (0.99) | 2.33 (0.99) |

**5. Empirical results and applications.** The empirical part of the paper includes in its first part a discussion on the computation of the critical values and the validity of the normal approximation provided in Remark 2.1 in the case of finite sample sizes (Section 5.1). In a second part, we assess the finite sample properties of the test procedures for several multivariate time series models (Section 5.2). Finally, we address the applicability by investigating volatilities and cross-volatilities of stock data for 12 companies (Section 5.3).

5.1. *Computation of critical values.* The asymptotic critical values of the test statistics $\Lambda_n$ and $\Omega_n$ used in Theorem 2.1 can be computed either using the precise distribution of the limits $\Lambda(\mathfrak{d})$ and $\Omega(\mathfrak{d})$, respectively, or—provided $\mathfrak{d}$ is sufficiently large—the normal approximation given in Remark 2.1. Critical values for $\mathfrak{d} = 1, \ldots, 5$ have already been tabulated in Kiefer [36]. Here, we provide critical values for larger values of $\mathfrak{d}$ that seem to be more realistic in today's financial applications. The calculations are based on the limit distributions of $\Lambda(\mathfrak{d})$ and $\Omega(\mathfrak{d})$ which have been derived in [36].

Let $q_{1-\alpha}^{\Omega}(\mathfrak{d})$ be the $1 - \alpha$ quantile of $\bar{\Omega}(\mathfrak{d})$, the standardized version of $\Omega(\mathfrak{d})$, and let furthermore $q_{1-\alpha}^{\Omega}(\infty) = z_{1-\alpha}$ denote the $1 - \alpha$ quantile of the standard normal distribution. In Table 1, we have recorded the critical values for three significance levels, namely $\alpha = 0.10$, 0.05 and 0.01, and various values of $\mathfrak{d}$ ranging from 10 to 500. Under each critical value, we report in brackets also the values $P(\bar{\Omega}(\mathfrak{d}) \leq z_{1-\alpha})$. Let $q_{1-\alpha}^{\Lambda}(\mathfrak{d})$ be the $1 - \alpha$ quantile of $\bar{\Lambda}(\mathfrak{d})$, the standardized version of $\Lambda(\mathfrak{d})$. The corresponding critical values for the same choices of $\alpha$ and $\mathfrak{d}$ have been collected in Table 2.

It can be seen from Table 1 that the normal approximation ($\mathfrak{d} = \infty$) works very reasonable in the case of $\bar{\Omega}(\mathfrak{d})$. Even for those quantiles $q_{1-\alpha}^{\Omega}$ that are the farthest from the corresponding standard normal quantile $z_{1-\alpha}$, it still holds that $P(\bar{\Omega}(\mathfrak{d}) \leq z_{1-\alpha}) \approx 1 - \alpha$. This, however, does not hold true for the second statistic $\bar{\Lambda}(\mathfrak{d})$. Here, the normal approximation ($\mathfrak{d} = \infty$) works considerably worse and it cannot even be recommended for $\mathfrak{d} = 500$. The exact values should consequently be used instead.

5.2. *A simulation study.* In this subsection, we report the results of a small simulation study that includes multivariate AR(1), dynamic factor model and exponential GARCH specifications for the process $(\mathbf{y}_j : j \in \mathbb{Z})$. In all cases, we use the statistics $\Omega_n$. Due to the weak dependence displayed by the data generating processes, we choose to work with the Bartlett estimator as a proxy for the asymptotic covariance matrix $\Sigma$ with window length $q = q(n) = \log_{10} n$ (see [4]). It is shown in Chapter 11 of [16] that the Bartlett estimator satisfies condition (2.6).

If $(\mathbf{y}_j : j \in \mathbb{Z})$ follows a parametric model, then $\Sigma$ is a function of the model parameters. These, in turn, yield a natural plug-in estimator $\hat{\Sigma}_n^*$ based on the relevant parameter estimators. This method, however, fails to work if the model is not correctly specified. Then, the plug-in estimator will converge to a limit $\Sigma^* \neq \Sigma$, resulting in a potential power loss of the test statistics. Moreover, estimation procedures and goodness of fit tests for the nonlinear models discussed in Section 4 still often lack accuracy and a nonparametric approach seems—at this point—more advisable.

AUTOREGRESSIVE PROCESSES. We used as data generating process the four-dimensional AR(1) time series given by the equations

$$\mathbf{y}_j = A\mathbf{y}_{j-1} + \boldsymbol{\varepsilon}_j, \qquad j = 1, \ldots, n,$$

TABLE 2
*Critical values for $\bar{\Lambda}(\mathfrak{d})$ for the three significance levels* 0.10, 0.05 *and* 0.01. *For each critical value, we report $P(\bar{\Lambda}(\mathfrak{d}) \leq z_{1-\alpha})$ in brackets. The last column, $\mathfrak{d} = \infty$, corresponds to the standard normal case as stated in Remark 2.1*

| $\mathfrak{d}$ | 10 | 15 | 20 | 50 | 100 | 200 | 500 | $\infty$ |
|---|---|---|---|---|---|---|---|---|
| $q_{0.90}^{\Lambda}(\mathfrak{d})$ | 2.64 (0.56) | 2.53 (0.58) | 2.46 (0.59) | 2.27 (0.63) | 2.16 (0.66) | 2.06 (0.69) | 1.96 (0.72) | 1.28 (0.90) |
| $q_{0.95}^{\Lambda}(\mathfrak{d})$ | 3.17 (0.69) | 3.02 (0.71) | 2.92 (0.72) | 2.69 (0.76) | 2.55 (0.78) | 2.44 (0.81) | 2.33 (0.83) | 1.64 (0.95) |
| $q_{0.99}^{\Lambda}(\mathfrak{d})$ | 4.28 (0.85) | 4.04 (0.87) | 3.89 (0.88) | 3.53 (0.91) | 3.33 (0.93) | 3.18 (0.94) | 3.04 (0.95) | 2.33 (0.99) |

where $\mathbf{y}_0$ has been obtained after a burn-in phase of 500 iterations. The sample sizes under consideration are $n = 200$, 500, 800 and 1000. The innovations $\boldsymbol{\varepsilon}_1, \ldots, \boldsymbol{\varepsilon}_n$ have been specified as independent 4-variate normal variates with independent components. To assess the empirical level of the test, we have chosen $A = (0.1)I$, where $I$ denotes the $4 \times 4$ identity matrix. Each component of the process is therefore a univariate AR(1) time series which is independent of the others.

To assess the power, we introduce several alternatives. These are specified by the form of the parameter matrix before and after the change-point $k^* = n/2$. In particular, we work here with $A = (0.1)I$ prior to $k^*$ and with $A^* = (0.1)I + \delta \mathbf{1}$ after $k^*$, where $\mathbf{1}$ denotes the $4 \times 4$ matrix for which all entries are equal to 1, thus introducing correlation between the coordinates after $k^*$. The values of $\delta$ vary between 0.1 and 0.6. Using Theorem 2.1, the critical values have been computed using the exact distribution of the limit random variable $\Omega(\mathfrak{d})$ given in (2.9) and discussed in the previous subsection. All results are based on $N = 1000$ repetitions.

The results are summarized in Table 3 for the three levels $\alpha = 0.10$, 0.05 and 0.01. The row with $\delta = 0$ corresponds to the empirical levels of the test. It can be seen that the procedure is conservative in finite samples as the empirical levels stay below the asymptotic critical values. Empirical and asymptotic level are close for the sample size $n = 1000$. On the other hand, small changes are found with greater difficulty if the sample size is small ($n = 200$). The power for bigger changes and/or larger sample sizes is, however, very reasonable.

FACTOR MODELS. As second data generating process we used a factor model in which the four-dimensional $\mathbf{y}_j$ are given by the equations

$$\mathbf{y}_j = L\mathbf{z}_j + \boldsymbol{\xi}_j, \qquad j = 1, \ldots, n,$$

TABLE 3
*The rejection levels of the statistic $\Omega_n$ for the various AR(1) processes of Section 5.2. The row specified by $\delta = 0$ contains the empirical levels*

| $\alpha$ | | 0.10 | | | | 0.05 | | | | 0.01 | | |
|---|---|---|---|---|---|---|---|---|---|---|---|---|
| $\delta/n$ | 200 | 500 | 800 | 1000 | 200 | 500 | 800 | 1000 | 200 | 500 | 800 | 1000 |
| 0.0 | 0.07 | 0.09 | 0.10 | 0.10 | 0.03 | 0.04 | 0.04 | 0.04 | 0.00 | 0.01 | 0.01 | 0.01 |
| 0.2 | 0.22 | 0.42 | 0.67 | 0.78 | 0.11 | 0.29 | 0.53 | 0.66 | 0.03 | 0.11 | 0.25 | 0.36 |
| 0.3 | 0.40 | 0.76 | 0.94 | 0.98 | 0.25 | 0.61 | 0.89 | 0.96 | 0.07 | 0.31 | 0.72 | 0.86 |
| 0.4 | 0.71 | 0.96 | 1.00 | 1.00 | 0.57 | 0.91 | 1.00 | 1.00 | 0.26 | 0.74 | 0.98 | 1.00 |
| 0.5 | 0.94 | 1.00 | 1.00 | 1.00 | 0.90 | 1.00 | 1.00 | 1.00 | 0.73 | 0.97 | 1.00 | 1.00 |
| 0.6 | 1.00 | 1.00 | 1.00 | 1.00 | 0.99 | 1.00 | 1.00 | 1.00 | 0.98 | 1.00 | 1.00 | 1.00 |

where $L$ is a $4 \times 2$ matrix, $\boldsymbol{\xi}_1, \ldots, \boldsymbol{\xi}_n$ are independent zero mean 4-variate normal variates with unit covariance matrix. The 2-variate factors $\mathbf{z}_1, \ldots, \mathbf{z}_n$ are given by the CCC–GARCH equations (4.3) and (4.4) with specifications

$$\boldsymbol{\omega} = \begin{pmatrix} 1.0 \\ 1.0 \end{pmatrix}, \qquad \boldsymbol{\alpha}_1 = \begin{pmatrix} 0.3 \\ 0.3 \end{pmatrix}, \qquad \boldsymbol{\beta}_1 = \begin{pmatrix} 0.3 \\ 0.3 \end{pmatrix}$$

and bivariate normal $\boldsymbol{\varepsilon}_j$ with identity covariance matrix. To assess the empirical level, we chose

$$L = \begin{pmatrix} 1 & 0 \\ 1 & 0 \\ 0 & 1 \\ 0 & 1 \end{pmatrix}$$

under $H_0$, while for the power considerations we picked $L^* = \delta L$ under $H_A$. The choices for the number of repetitions $N$, the sample sizes $n$ and the significance level $\alpha$ are the same as for the AR(1) processes above, while for the change sizes we used values of $\delta$ ranging from 1.1 to 2.0. The results are summarized in Table 4. The empirical levels appear to be more volatile than in the AR(1) case but quite satisfactory for the sample size $n = 1000$. The findings for the power of the test are similar to what has been observed for the AR(1) time series.

EXPONENTIAL GARCH MODELS. The final data generating process is a four-dimensional exponential GARCH model given by the (4.7) and (4.10) with specifications $p = q = 1$, $A_1 = 0.1 \operatorname{diag}(\mathbb{I}_{(10)})$, $B_1^T = (0.1 \operatorname{diag}(\mathbb{I}_{(4)}), \mathbf{0})$ and $F = 0$, where $\mathbb{I}_{(d)}$ denotes the $d$-dimensional vector whose entries are all equal to 1 and $\mathbf{0}$ the $4 \times 4$ zero matrix. Under $H_0$, we worked with $C = 0.2 \operatorname{diag}(\mathbb{I}_{(4)})$, while under $H_A$ we used $C$ prior to $k^*$ and with $C^* = \delta C$ thereafter, where $\delta$ ranges from 1 to 3.5. The errors $(\boldsymbol{\varepsilon}_j : j = 1, \ldots, n)$ are standard normal and the choices for $n$, $k^*$ and $\alpha$ are the same as before. The results are summarized in Table 5. They are close to the findings in the AR(1) and factor model cases.

TABLE 4
*The rejection levels of the statistic $\Omega_n$ for the factor model of Section 5.2. The row specified by $\delta = 1$ contains the empirical levels*

| $\alpha$ | **0.10** | | | | **0.05** | | | | **0.01** | | | |
|---|---|---|---|---|---|---|---|---|---|---|---|---|
| $\delta/n$ | **200** | **500** | **800** | **1000** | **200** | **500** | **800** | **1000** | **200** | **500** | **800** | **1000** |
| 1.0 | 0.15 | 0.12 | 0.13 | 0.15 | 0.06 | 0.05 | 0.07 | 0.08 | 0.01 | 0.01 | 0.01 | 0.01 |
| 1.1 | 0.16 | 0.19 | 0.26 | 0.31 | 0.09 | 0.09 | 0.16 | 0.21 | 0.02 | 0.02 | 0.05 | 0.06 |
| 1.2 | 0.25 | 0.39 | 0.58 | 0.66 | 0.13 | 0.27 | 0.45 | 0.55 | 0.03 | 0.12 | 0.23 | 0.30 |
| 1.4 | 0.52 | 0.85 | 0.96 | 0.99 | 0.37 | 0.75 | 0.93 | 0.98 | 0.13 | 0.51 | 0.83 | 0.92 |
| 1.6 | 0.76 | 0.98 | 1.00 | 1.00 | 0.63 | 0.96 | 1.00 | 1.00 | 0.35 | 0.86 | 0.99 | 1.00 |
| 1.8 | 0.91 | 1.00 | 1.00 | 1.00 | 0.82 | 1.00 | 1.00 | 1.00 | 0.56 | 0.99 | 1.00 | 1.00 |
| 2.0 | 0.96 | 1.00 | 1.00 | 1.00 | 0.93 | 1.00 | 1.00 | 1.00 | 0.74 | 0.99 | 1.00 | 1.00 |

TABLE 5
*The rejection levels of the statistic $\Omega_n$ for the exponential GARCH model of Section 5.2. The row specified by $\delta = 1$ contains the empirical levels*

| $\alpha$ | 0.10 | | | | 0.05 | | | | 0.01 | | | |
|---|---|---|---|---|---|---|---|---|---|---|---|---|
| $\delta/n$ | 200 | 500 | 800 | 1000 | 200 | 500 | 800 | 1000 | 200 | 500 | 800 | 1000 |
| 1.0 | 0.07 | 0.07 | 0.07 | 0.08 | 0.04 | 0.03 | 0.04 | 0.04 | 0.00 | 0.01 | 0.00 | 0.01 |
| 1.5 | 0.10 | 0.19 | 0.26 | 0.30 | 0.04 | 0.09 | 0.16 | 0.18 | 0.00 | 0.01 | 0.04 | 0.06 |
| 2.0 | 0.18 | 0.47 | 0.75 | 0.83 | 0.08 | 0.31 | 0.61 | 0.74 | 0.01 | 0.11 | 0.34 | 0.49 |
| 2.5 | 0.40 | 0.83 | 0.98 | 0.99 | 0.25 | 0.74 | 0.96 | 0.98 | 0.07 | 0.47 | 0.85 | 0.95 |
| 3.0 | 0.61 | 0.98 | 1.00 | 1.00 | 0.43 | 0.96 | 1.00 | 1.00 | 0.17 | 0.86 | 1.00 | 1.00 |
| 3.5 | 0.83 | 1.00 | 1.00 | 1.00 | 0.70 | 1.00 | 1.00 | 1.00 | 0.38 | 0.98 | 1.00 | 1.00 |

ESTIMATING THE BREAK LOCATION. To conclude the simulation study, we briefly investigated also how well the relative change location $\theta$ is approximated by the estimator $\hat{\theta}_n$ of (2.12) for the AR(1) processes, the factor and exponential GARCH models. In Table 6, we have reported mean, median and standard deviation for the various choices of $\delta$ and $n$ specified in the foregoing paragraphs and for the break point $k^* = n/2$, thus leading to $\theta = 0.5$. The significance level was set to $\alpha = 0.1$. The table shows that for all three processes there is a positive bias which diminishes as the sample size or the change magnitude $\delta$ increase. Generally, using the median leads to a superior approximation of the true $\theta$ than utilizing the mean, thereby indicating the existence of unusually large deviations (outliers) in several of the simulation runs. The estimator $\hat{\theta}_n$ works satisfactory in the case $n = 1000$. This shows that the asymptotic result in display (2.13) applies reasonably well here, although future research should pursue a more detailed in-depth analysis of the performance of $\hat{\theta}_n$.

5.3. *An application to financial data.* We study in our application the log returns of the adjusted closing stock prices between July 19, 1993, and March 19, 2009, of the 12 companies listed in Table 7, yielding the 12-dimensional observations $\mathbf{y}_1, \ldots, \mathbf{y}_{3941}$. There are four companies in the airline sector, four in the automotive sector and four in the energy sector. More precisely, if $p_{j,\ell}$ denotes the price of stock $\ell$ at time $j$, then we work with the centered log-returns

$$y_{j,\ell} = \log\left(\frac{p_{j+1,\ell}}{p_{j,\ell}}\right) - \frac{1}{3941}\sum_{i=1}^{3941}\log\left(\frac{p_{i+1,\ell}}{p_{i,\ell}}\right),$$

TABLE 6
*The estimation of the change-point fraction $\theta$ for the AR(1) processes (upper panel), the factor GARCH model (middle panel) and the exponential GARCH model (lower panel)*

| $n$ | 200 | | | 500 | | | 800 | | | 1000 | | |
|-----|------|------|------|------|------|------|------|------|------|------|------|------|
| $\delta$ | mean | med | sd | mean | med | sd | mean | med | sd | mean | med | sd |
| 0.2 | 0.55 | 0.55 | 0.11 | 0.56 | 0.55 | 0.09 | 0.54 | 0.53 | 0.08 | 0.54 | 0.52 | 0.07 |
| 0.3 | 0.56 | 0.55 | 0.10 | 0.55 | 0.53 | 0.08 | 0.54 | 0.52 | 0.06 | 0.53 | 0.52 | 0.05 |
| 0.4 | 0.59 | 0.57 | 0.09 | 0.55 | 0.53 | 0.06 | 0.53 | 0.52 | 0.05 | 0.53 | 0.51 | 0.04 |
| 0.5 | 0.58 | 0.57 | 0.08 | 0.54 | 0.52 | 0.05 | 0.53 | 0.52 | 0.05 | 0.53 | 0.51 | 0.04 |
| 1.2 | 0.56 | 0.56 | 0.11 | 0.54 | 0.53 | 0.09 | 0.54 | 0.52 | 0.09 | 0.53 | 0.52 | 0.08 |
| 1.4 | 0.56 | 0.55 | 0.09 | 0.53 | 0.52 | 0.07 | 0.52 | 0.51 | 0.05 | 0.52 | 0.51 | 0.04 |
| 1.6 | 0.55 | 0.53 | 0.07 | 0.52 | 0.51 | 0.05 | 0.51 | 0.51 | 0.03 | 0.51 | 0.51 | 0.03 |
| 1.8 | 0.53 | 0.52 | 0.06 | 0.52 | 0.51 | 0.03 | 0.51 | 0.51 | 0.02 | 0.51 | 0.50 | 0.02 |
| 2.0 | 0.53 | 0.52 | 0.05 | 0.51 | 0.50 | 0.03 | 0.51 | 0.50 | 0.02 | 0.51 | 0.50 | 0.02 |
| 1.5 | 0.52 | 0.54 | 0.12 | 0.54 | 0.53 | 0.09 | 0.52 | 0.52 | 0.09 | 0.53 | 0.52 | 0.09 |
| 2.0 | 0.54 | 0.53 | 0.10 | 0.53 | 0.52 | 0.07 | 0.52 | 0.51 | 0.07 | 0.52 | 0.51 | 0.06 |
| 2.5 | 0.54 | 0.53 | 0.09 | 0.52 | 0.51 | 0.04 | 0.51 | 0.50 | 0.04 | 0.51 | 0.50 | 0.03 |
| 3.0 | 0.53 | 0.52 | 0.07 | 0.51 | 0.50 | 0.02 | 0.51 | 0.50 | 0.02 | 0.51 | 0.50 | 0.02 |
| 3.5 | 0.53 | 0.51 | 0.06 | 0.51 | 0.50 | 0.02 | 0.51 | 0.50 | 0.02 | 0.50 | 0.50 | 0.01 |

$$j = 1, \ldots, 3941; \ell = 1, \ldots, 12.$$

To illustrate the time-varying nature of the volatilities, we have computed the rolling volatility and cross-volatility estimators

$$\hat{\gamma}_j(k, \ell) = \frac{1}{100} \sum_{i=j-100+1}^{j} y_{i,k} y_{i,\ell}, \qquad j = 101, \ldots, 3941; k, \ell = 1, \ldots, 12.$$

For each of the three sectors, the corresponding estimators are given in Figure 1. They appear to be time-dependent. To further assess this conjecture, we computed the test statistic value $\Omega_{3941} = 60.07$. Since $d = 12$, we have $\mathfrak{d} = 78$. This leads to the critical value 12.00 at the 95% level. Therefore, the hypothesis of no change in the volatility and cross-volatility structure is rejected.

To illustrate the potential to detect multiple departures of the stability hypothesis, we have applied a binary segmentation procedure of the original test. Each time we reject $H_0$, we re-apply it to the subsamples obtained from splitting the data into two pieces at the observation determined by (2.12), that is, $\hat{k}^* = \lfloor \hat{\theta}_n n \rfloor$. The corresponding findings are reported in Table 8. A number of the detected changes can be readily associated with major historical events. For example, the change found in the first round of the segmentation procedure along with the other estimated changes in 2001 are linked to the bursting of the dot-com bubble and the September 11 attacks, while the break dates in 1997 and 1998 are connected to the Asian financial crisis, and the collapse of the hedge fund Long-Term Capital Management and the Russian financial crisis, respectively. The detected breaks in 2007 and 2008 can be related to the collapse of the housing market in the United States and several European countries. Finally, the latest significant break, which was found at lag 3809, signifies the beginning of the recent financial crisis. The detected change-point time is September 9, 2008, thus predating the collapse of the investment bank Lehman Brothers by only three trading days.

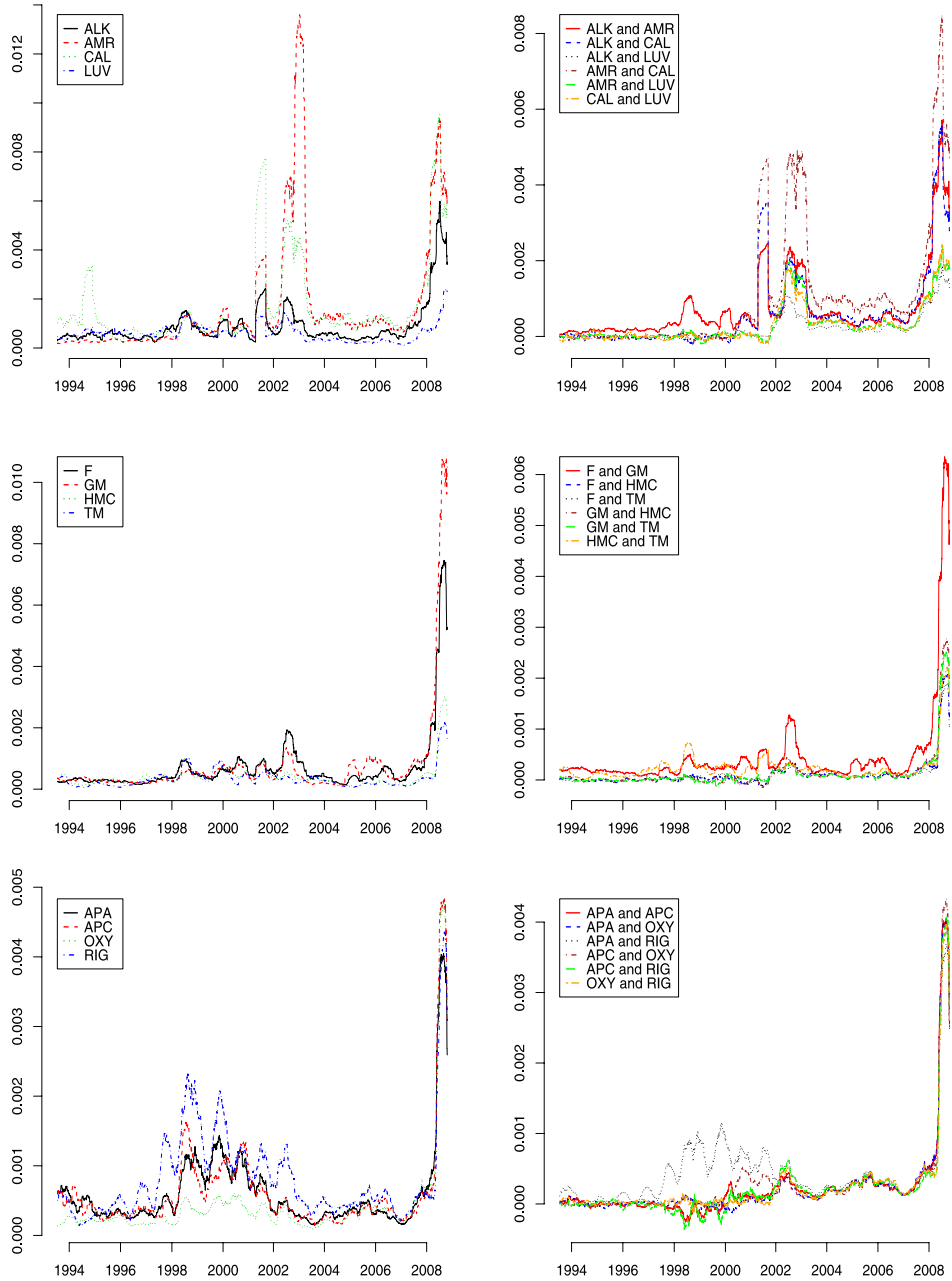

FIG. 1. *The volatilities (left) and cross-volatilities (right) of the airline sector (upper panel), the automotive sector (middle panel), and the energy sector (lower panel). The company abbreviations used in the plots are detailed in Table 7.*

TABLE 7
*The 12 companies whose stock values are studied in Section 5.3*

|  | Symbol | Name | Sector |
|---|---|---|---|
| 1 | ALK | Alaska Air Group, Inc. | Airline |
| 2 | AMR | AMR Corp. (American Airlines) | Airline |
| 3 | CAL | Continental Airlines, Inc. | Airline |
| 4 | LUV | Southwest Airlines Co. | Airline |
| 5 | F | Ford Motor Co. | Automotive |
| 6 | GM | General Motors Corp. | Automotive |
| 7 | HMC | Honda Motor Co. | Automotive |
| 8 | TM | Toyota Motor Corp. | Automotive |
| 9 | APA | Apache Corp. | Energy |
| 10 | APC | Anadarko Petroleum Corp. | Energy |
| 11 | OXY | Occidental Petroleum Corp. | Energy |
| 12 | RIG | Transocean Inc. | Energy |

TABLE 8
*The segmentation procedure performed on the entire data set. The estimated change-points are listed with their corresponding test statistic value, the round in which the change was found and whether or not it was significant at the 95% error level*

| $k^*$ | $k^*$ (Date) | $\Omega_n$ | Round | Significant |
|---|---|---|---|---|
| 305 | 1994-09-30 | 11.60 | 4 | No |
| 647 | 1996-02-07 | 19.80 | 3 | Yes |
| 1021 | 1997-07-31 | 14.78 | 4 | Yes |
| 1259 | 1998-07-13 | 34.05 | 2 | Yes |
| 1400 | 1999-02-02 | 13.60 | 4 | Yes |
| 1770 | 2000-07-25 | 18.65 | 3 | Yes |
| 1886 | 2001-01-11 | 12.07 | 4 | Yes |
| 2101 | 2001-11-23 | 60.07 | 1 | Yes |
| 2358 | 2002-12-03 | 14.66 | 4 | Yes |
| 2728 | 2004-05-24 | 27.90 | 3 | Yes |
| 3175 | 2006-03-03 | 14.70 | 4 | Yes |
| 3589 | 2007-10-24 | 28.95 | 2 | Yes |
| 3710 | 2008-04-18 | 10.11 | 4 | No |
| 3809 | 2008-09-09 | 12.13 | 3 | Yes |
| 3871 | 2008-12-05 | 7.47 | 4 | No |

## APPENDIX A: A MULTIVARIATE FUNCTIONAL CENTRAL LIMIT THEOREM

In this section, we derive a general multivariate functional central limit theorem (FCLT) and state moreover a derivative that is tailor-made for the examples discussed in Section 4. The results extend the corresponding univariate theory provided, for example, in Billingsley's monograph [10].

Recall from Assumption 2.1 that $(\boldsymbol{\varepsilon}_j : j \in \mathbb{Z})$ is a sequence of independent, identically distributed random variables with values in $\mathbb{R}^{d'}$ and that we further suppose that the $d$-dimensional sequence $(\mathbf{y}_j : j \in \mathbb{Z})$ is given by the relations (2.3). To prove a FCLT for this sequence, we shall utilize the sequence of $m$-dependent random vectors $(\mathbf{y}_j^{(m)} : j \in \mathbb{Z})$ given in (2.4), which is sufficiently close to $(\mathbf{y}_j : j \in \mathbb{Z})$ in the sense of condition (A.1) below. Recall finally that a stochastic process $(\mathbf{W}_\Gamma(t) : t \in [0, 1])$ is called a $d$-dimensional Brownian motion with covariance matrix $\Gamma \in \mathbb{R}^{d \times d}$ if it is Gaussian with mean $\mathbf{0}$ and covariance function $\mathrm{Cov}(\mathbf{W}_\Gamma(s), \mathbf{W}_\Gamma(t)) = \min\{s, t\}\Gamma$.

THEOREM A.1. *Assume that the $d$-dimensional random process $(\mathbf{y}_j : j \in \mathbb{Z})$ specified by (2.3) satisfies $E[\mathbf{y}_j] = 0$ and $E[|\mathbf{y}_j|^2] < \infty$. Suppose further that, for any $m \geq 1$, the vectors $\mathbf{y}_0^{(m)}$ can be defined such that*

$$(A.1) \qquad \sum_{m \geq 1} \|\mathbf{y}_0 - \mathbf{y}_0^{(m)}\|_2 < \infty.$$

*Then the series $\Gamma = \sum_{j \in \mathbb{Z}} \mathrm{Cov}(\mathbf{y}_0, \mathbf{y}_j)$ converges (coordinatewise) absolutely and*

$$\frac{1}{\sqrt{n}} \sum_{j=1}^{[nt]} \mathbf{y}_j \overset{D_d[0,1]}{\longrightarrow} \mathbf{W}_\Gamma(t) \qquad (n \to \infty),$$

*where the convergence takes place in the $d$-dimensional Skorohod space $D_d[0, 1]$.*

The proof of Theorem A.1 needs the following small auxiliary result.

LEMMA A.1. *Let $(\mathbf{X}_n : n \geq 1)$ and $(\mathbf{Y}_n : n \geq 1)$ be sequences of random vectors in $\mathbb{R}^d$ such that $\mathbf{X}_n \overset{\mathcal{D}}{\longrightarrow} \mathbf{X}$ and $\mathbf{Y}_n \overset{\mathcal{D}}{\longrightarrow} \mathbf{Y}$ as $n \to \infty$, where the limit random vectors $\mathbf{X}$ and $\mathbf{Y}$ are independent. Suppose further that there are sequences $(\mathbf{X}_n^* : n \geq 1)$ and $(\mathbf{Y}_n^* : n \geq 1)$ such that (a) $\mathbf{X}_n^*$ and $\mathbf{Y}_n^*$ are independent for all $n \geq 1$, and (b) $|\mathbf{X}_n^* - \mathbf{X}_n| \overset{\mathcal{P}}{\longrightarrow} 0$ and $|\mathbf{Y}_n^* - \mathbf{Y}_n| \overset{\mathcal{P}}{\longrightarrow} 0$ as $n \to \infty$. Then it holds that $(\mathbf{X}_n, \mathbf{Y}_n) \overset{\mathcal{D}}{\longrightarrow} (\mathbf{X}, \mathbf{Y})$ as $n \to \infty$.*

PROOF. By assumption, it holds that $\mathbf{X}_n^* \xrightarrow{\mathcal{D}} \mathbf{X}$ and $\mathbf{Y}_n^* \xrightarrow{\mathcal{D}} \mathbf{Y}$. Since the sequences $(\mathbf{X}_n^* : n \geq 1)$ and $(\mathbf{Y}_n^* : n \geq 1)$ are independent, we obtain $(\mathbf{X}_n^*, \mathbf{Y}_n^*) \xrightarrow{\mathcal{D}} (\mathbf{X}, \mathbf{Y})$. Viewing $(\mathbf{X}_n^*, \mathbf{Y}_n^*)$ and $(\mathbf{X}, \mathbf{Y})$ as random vectors in $\mathbb{R}^{2d}$ implies $|(\mathbf{X}_n^*, \mathbf{Y}_n^*) - (\mathbf{X}_n, \mathbf{Y}_n)| \xrightarrow{\mathcal{D}} 0$ and consequently the assertion. □

PROOF OF THEOREM A.1.

*Step* 1. For $n \geq 1$ and $t \in [0,1]$, let $\mathcal{V}_n(t) = n^{-1/2} \sum_{j=1}^{[nt]} \mathbf{y}_j$. It will be shown that the finite-dimensional distributions of the partial sum process $(\mathcal{V}_n(t) : t \in [0,1])$ converge to the finite-dimensional distributions corresponding to the $d$-dimensional Brownian motion $(\mathbf{W}_\Sigma(t) : t \in [0,1])$. Let $K \geq 1$ and $0 \leq t_0 < t_1 < \cdots < t_K \leq 1$. Because the limit process has independent increments, proving the latter convergence is equivalent to showing that

$$
\begin{aligned}
(A.2) \qquad & (\mathcal{V}_n(t_{\ell-1}, t_\ell) : 1 \leq \ell \leq K) \\
& \xrightarrow{\mathcal{D}} (\mathbf{W}_\Gamma(t_\ell) - \mathbf{W}_\Gamma(t_{\ell-1}) : 1 \leq \ell \leq K) \qquad (n \to \infty),
\end{aligned}
$$

where $\mathcal{V}_n(s,t) = n^{-1/2} \sum_{j=[ns]+1}^{[nt]} \mathbf{y}_j$ for $0 \leq s < t \leq 1$. To do so, we focus in the first part of the proof on one element of the random vectors in (A.2). More precisely, we shall establish with the Cramér–Wold device that, for any $0 \leq s < t \leq 1$,

$$
(A.3) \qquad \mathcal{V}_n(s,t) \xrightarrow{\mathcal{D}} \mathbf{W}_\Gamma(t) - \mathbf{W}_\Gamma(s), \qquad (n \to \infty).
$$

For $\mathbf{a} = (a_1, \ldots, a_d)^T \in \mathbb{R}^d$, let the univariate sequences $(z_j : j \in \mathbb{Z})$ and $(z_j^{(m)} : j \in \mathbb{Z})$ be defined by letting

$$
z_j = \mathbf{a}^T \mathbf{y}_j \quad \text{and} \quad z_j^{(m)} = \mathbf{a}^T \mathbf{y}_j^{(m)},
$$

respectively. Let $\mathbf{y}_0 = (y_{0,1}, \ldots, y_{0,d})^T$ and $\mathbf{y}_0^{(m)} = (y_{0,1}^{(m)}, \ldots, y_{0,d}^{(m)})^T$ and observe that $|y_{0,i} - y_{0,i}^{(m)}| \leq |\mathbf{y}_0 - \mathbf{y}_0^{(m)}|$ for all $1 \leq i \leq d$. Therefore,

$$
\begin{aligned}
\sum_{m \geq 1} \|z_0 - z_0^{(m)}\|_2 &\leq \sum_{m \geq 1} \max_{1 \leq i \leq d} |a_i| \|\mathbf{y}_0 - \mathbf{y}_0^{(m)}\|_2 \\
&\leq \left( \sum_{i=1}^d |a_i| \right) \left( \sum_{m \geq 1} \|\mathbf{y}_0 - \mathbf{y}_0^{(m)}\|_2 \right) < \infty
\end{aligned}
$$

by assumption (A.1). Theorem 21.1 in [10] consequently yields that

$$
\frac{1}{\sqrt{n}} \sum_{j=[ns]+1}^{[nt]} z_j \xrightarrow{\mathcal{D}} N(0, \sigma_z^2(t-s)) \qquad (n \to \infty),
$$

where $\sigma_z^2 = \sum_{j \in \mathbb{Z}} E[z_0 z_j] = \sum_{j \in \mathbb{Z}} \mathbf{a}^T \operatorname{Cov}(\mathbf{y}_0, \mathbf{y}_j) \mathbf{a} < \infty$, which implies in particular that $\Gamma$ converges coordinatewise absolutely. On the other hand, it holds that

$$\mathbf{a}^T(\mathbf{W}_\Gamma(t) - \mathbf{W}_\Gamma(s)) \overset{\mathcal{D}}{=} N(0, \sigma_z^2(t-s)),$$

so that (A.3) is proved. However, since the partial sums $\mathcal{V}_n(t_{\ell-1}, t_\ell)$, $1 \leq \ell \leq K$, are dependent, (A.3) does not yet imply (A.2). In the second part of the proof, we shall therefore construct independent random variables to replace them. For $0 \leq s < t \leq 1$, let $\mathcal{V}_n^*(s,t) = n^{-1/2} \sum_{j=[ns]+1}^{[nt]} \mathbf{y}_j^{([ns]+1-j)}$, where $\mathbf{y}_j^{(0)} = \mathbf{0}$. Then $\mathcal{V}_n^*(s,t)$ and $\mathcal{V}_n(0,s)$ are independent. Moreover, for any $\delta > 0$, it holds that

$$P(|\mathcal{V}_n(s,t) - \mathcal{V}_n^*(s,t)| > \delta) \leq \frac{1}{\delta^2} E[|\mathcal{V}_n(s,t) - \mathcal{V}_n^*(s,t)|^2]$$

$$\leq \frac{1}{\delta^2 n} \left( \sum_{m=0}^{[nt]-[ns]} \|\mathbf{y}_{[ns]+m} - \mathbf{y}_{[ns]+m}^{(m)}\|_2 \right)^2$$

$$\leq \frac{1}{\delta^2 n} \left( \sum_{m=0}^{\infty} \|\mathbf{y}_0 - \mathbf{y}_0^{(m)}\|_2 \right)^2 \to 0$$

as $n \to \infty$. An application of Lemma A.1 shows that the vectors $(\mathcal{V}_n(0,s), \mathcal{V}_n(s,t))$ converge in distribution to $(\mathbf{W}_\Gamma(s), \mathbf{W}_\Gamma(t) - \mathbf{W}_\Gamma(s))$ as $n \to \infty$. This proves (A.2) for $K = 2$. The same arguments apply also for arbitrary $K$.

*Step* 2. It remains to verify that the partial sum process $(\mathcal{V}_n(t): t \in [0,1])$ is tight. Since it suffices to show the tightness coordinatewise, a reference to Theorem 21.1 in [10] completes the proof of Theorem A.1.  $\square$

The main results in Section 2 cover the case of detecting breaks in the covariance structure of multivariate financial time series. Therefore, at least in the context of the present paper, the following derivative of Theorem A.1 is of greater importance.

THEOREM A.2.  *Suppose that the assumptions of Theorem A.1 hold true and that* $E[|\mathbf{y}_j|^4] < \infty$. *If* (2.5) *is satisfied, then the series* $\Sigma = \sum_{j \in \mathbb{Z}} \operatorname{Cov}(\operatorname{vech}[\mathbf{y}_0 \mathbf{y}_0^T], \operatorname{vech}[\mathbf{y}_j \mathbf{y}_j^T])$ *converges (coordinatewise) absolutely and*

$$\frac{1}{\sqrt{n}} \sum_{j=1}^{[nt]} (\operatorname{vech}[\mathbf{y}_j \mathbf{y}_j^T - E[\mathbf{y}_j \mathbf{y}_j^T]]) \overset{D_{\mathfrak{d}}[0,1]}{\longrightarrow} \mathbf{W}_\Sigma(t) \qquad (n \to \infty),$$

*where* $\mathfrak{d} = d(d+1)/2$ *and the convergence takes place in the* $\mathfrak{d}$-*dimensional Skorohod space* $D_{\mathfrak{d}}[0,1]$.

PROOF. First, we define the symmetric matrix $Y_0 = \mathbf{y}_0\mathbf{y}_0^T - E[\mathbf{y}_0\mathbf{y}_0^T]$. To prove the assertion of Theorem A.2, it suffices to show that the assumptions of Theorem A.1 are satisfied for

$$\mathbf{z}_0 = \mathrm{vech}[Y_0].$$

Note that the moment assumptions $E[\mathbf{z}_0] = 0$ and $E[|\mathbf{z}_0|^2] < \infty$ are clearly fulfilled, so that it remains to verify condition (A.1). Recall that $\mathbf{y}_0 = (y_{0,1}, \ldots, y_{0,d})^T$ and $\mathbf{y}_0^{(m)} = (y_{0,1}^{(m)}, \ldots, y_{0,d}^{(m)})^T$. The $(k, \ell)$th element of $Y_0$ is thus given by

$$Y_0(k, \ell) = y_{0,k}y_{0,\ell} - E[y_{0,k}y_{0,\ell}], \qquad 1 \le k, \ell \le d.$$

In a similar fashion, we define the matrix $Y_0^{(m)}$ corresponding to the approximating variable $\mathbf{z}_0^{(m)} = \mathrm{vech}[Y_0^{(m)}]$ by letting

$$Y_0^{(m)}(k, \ell) = y_{0,k}^{(m)} y_{0,\ell}^{(m)} - E[y_{0,k}y_{0,\ell}], \qquad 1 \le k, \ell \le d.$$

Observe next that, applying Cauchy–Schwarz and Minkowski, with some constant $c > 0$,

$$E[|y_{0,k}y_{0,\ell} - y_{0,k}^{(m)} y_{0,\ell}^{(m)}|^2]$$

$$\le E[(|y_{0,k} - y_{0,k}^{(m)}||y_{0,\ell}| + |y_{0,k}^{(m)}||y_{0,\ell} - y_{0,\ell}^{(m)}|)^2]$$

$$\le \max_{1 \le k, \ell \le d}\{\|y_{0,\ell}\|_4^2, \|y_{0,k}^{(m)}\|_4^2\}(\|y_{0,k} - y_{0,k}^{(m)}\|_4 + \|y_{0,\ell} - y_{0,\ell}^{(m)}\|_4)^2$$

$$\le 4 \max_{1 \le k, \ell \le d}\{\|y_{0,\ell}\|_4^2, (\|y_{0,k}\|_4 + \|y_{0,k} - y_{0,k}^{(m)}\|_4)^2\}\|\mathbf{y}_0 - \mathbf{y}_0^{(m)}\|_4^2$$

$$\le 4\left(\|\mathbf{y}_0\|_4 + \sum_{m=1}^{\infty}\|\mathbf{y}_0 - \mathbf{y}_0^{(m)}\|_4\right)^2\|\mathbf{y}_0 - \mathbf{y}_0^{(m)}\|_4^2$$

$$= c^2\|\mathbf{y}_0 - \mathbf{y}_0^{(m)}\|_4^2,$$

since $\|y_{0,k} - y_{0,k}^{(m)}\|_4 \le \|\mathbf{y}_0 - \mathbf{y}_0^{(m)}\|_4 \le \sum_{m=1}^{\infty}\|\mathbf{y}_0 - \mathbf{y}_0^{(m)}\|_4 < \infty$. Thus, we arrive at

$$\sum_{m \ge 1}\|\mathbf{z}_0 - \mathbf{z}_0^{(m)}\|_2 \le \sum_{m=1}^{\infty}\left(\sum_{k=1}^{d}\sum_{\ell=1}^{d} E[|y_{0,k}y_{0,\ell} - y_{0,k}^{(m)} y_{0,\ell}^{(m)}|^2]\right)^{1/2}$$

$$\le cd \sum_{m=1}^{\infty}\|\mathbf{y}_0 - \mathbf{y}_0^{(m)}\|_4 < \infty$$

and the proof is therefore complete. $\square$

## APPENDIX B: MATHEMATICAL PROOFS

**B.1. Proofs of Theorems 2.1, 2.2 and Remark 2.1.**

PROOF OF THEOREM 2.1.   Theorem A.2 yields that

$$\frac{1}{\sqrt{n}} \sum_{j=1}^{\lfloor nt \rfloor} (\text{vech}[\mathbf{y}_j \mathbf{y}_j^T] - E[\text{vech}[\mathbf{y}_j \mathbf{y}_j^T]]) \xrightarrow{D^{\mathfrak{d}}[0,1]} \mathbf{W}_\Sigma(t) \qquad (n \to \infty).$$

Hence,

$$\frac{1}{\sqrt{n}} \left( \sum_{j=1}^{\lfloor nt \rfloor} \text{vech}[\mathbf{y}_j \mathbf{y}_j^T] - t \sum_{j=1}^{n} \text{vech}[\mathbf{y}_j \mathbf{y}_j^T] \right) \xrightarrow{D^{\mathfrak{d}}[0,1]} \mathbf{B}_\Sigma(t) \qquad (n \to \infty),$$

where $(\mathbf{B}_\Sigma(t) : t \in [0,1])$ denotes a Gaussian process with mean function $E[\mathbf{B}_\Sigma(t)] = \mathbf{0}$ and covariance function $E[\mathbf{B}_\Sigma(t) \mathbf{B}_\Sigma^T(s)] = \Sigma(\min\{t,s\} - ts)$. Now Theorem 2.1 follows from the continuous mapping theorem and assumption (2.6).   □

PROOF OF REMARK 2.1.   We consider $\Lambda(\mathfrak{d})$ first and show in a first step that

$$(\text{B.1}) \qquad \Lambda_{\mathfrak{d}}^*(t) = \frac{1}{\sqrt{\mathfrak{d}}} \sum_{\ell=1}^{\mathfrak{d}} (B_\ell^2(t) - t[1-t]) \xrightarrow{C[0,1]} \Lambda^*(t) \qquad (\mathfrak{d} \to \infty),$$

where $\xrightarrow{C[0,1]}$ indicates weak convergence in the space of continuous functions defined on the interval $[0,1]$ and $(\Lambda^*(t) : t \in [0,1])$ denotes a continuous Gaussian process. Note that in the discourse of the proof, we shall not need the covariance structure of this process and we consequently do not compute it. We are going to apply the Cramér–Wold device. It follows from the independence of the Brownian bridge processes $(B_\ell(t) : t \in [0,1])$, $1 \le \ell \le \mathfrak{d}$, that, for all $t_i \in [0,1]$ and $\lambda_i \in \mathbb{R}$, $1 \le i \le N$,

$$\sum_{i=1}^{N} \lambda_i \Lambda_{\mathfrak{d}}^*(t_i) \xrightarrow{\mathcal{D}} \xi \qquad (\mathfrak{d} \to \infty),$$

where $\xi$ denotes a centered normal random variable. Next, we verify the validity of the moment condition (12.51) in Billingsley [10]. Let $\mu(t,s) = t(1-t) - s(1-s)$ and denote in the following by $c$ a universal constant which may vary from line to line. Observe that Rosenthal's inequality (see, e.g., page 59 of Petrov [44]) implies that

$$E[(\Lambda_{\mathfrak{d}}^*(t) - \Lambda_{\mathfrak{d}}^*(s))^4]$$
$$\le \frac{c}{\mathfrak{d}^2} \left[ \sum_{\ell=1}^{\mathfrak{d}} E[(B_\ell^2(t) - B_\ell^2(s) - \mu(t,s))^4] \right.$$

$$+ \left( \sum_{\ell=1}^{\mathfrak{d}} E[(B_\ell^2(t) - B_\ell^2(s) - \mu(t,s))^2] \right)^2 \right] \le c(t-s)^2,$$

since, by the Cauchy–Schwarz inequality and the fact that $E[(B_\ell(t) - B_\ell(s))^8] \le c(t-s)^4$, we have after routine computations that

$$E[(B_\ell^2(t) - B_\ell^2(s) - \mu(t,s))^4] \le c(t-s)^2,$$

$$E[(B_\ell^2(t) - B_\ell^2(s) - \mu(t,s))^2] \le c|t-s|.$$

Theorem 12.3 in [10] yields now the tightness and the continuity of the limit process $(\Lambda^*(t) : t \in [0,1])$.

Let $t_{\mathfrak{d}}$ denote the location of the (first) maximum of $B_1^2(t) + \cdots + B_{\mathfrak{d}}^2(t)$. Since the function $t \mapsto t(1-t)$ reaches its largest value for $t = 1/2$, the weak convergence in (B.1) implies that

(B.2) $$\left| t_{\mathfrak{d}} - \tfrac{1}{2} \right| = \mathcal{O}_P(\mathfrak{d}^{-1/4}) \qquad (\mathfrak{d} \to \infty).$$

We show now that

(B.3) $$\left| \sup_{0 \le t \le 1} \sum_{\ell=1}^{\mathfrak{d}} B_\ell^2(t) - \sum_{\ell=1}^{\mathfrak{d}} B_\ell^2\left( \frac{1}{2} \right) \right| = o_P(\mathfrak{d}^{1/2}) \qquad (\mathfrak{d} \to \infty).$$

Using (B.2), for any $\delta \in (0,1)$, there are a constant $c$ and a positive integer $\mathfrak{d}_0$ such that

$$P\left( \sup_{0 \le t \le 1} \sum_{\ell=1}^{\mathfrak{d}} B_\ell^2(t) = \sup_{|t-1/2| \le c/\mathfrak{d}^{1/4}} \sum_{\ell=1}^{\mathfrak{d}} B_\ell^2(t) \right) \ge 1 - \delta,$$

if $\mathfrak{d} \ge \mathfrak{d}_0$. Now

$$0 \le \sup_{|t-1/2| \le c/\mathfrak{d}^{1/4}} \left( \sum_{\ell=1}^{\mathfrak{d}} B_\ell^2(t) - \sum_{\ell=1}^{\mathfrak{d}} B_\ell^2\left( \frac{1}{2} \right) \right)$$

$$= \sup_{|t-1/2| \le c/\mathfrak{d}^{1/4}} \left( \sqrt{\mathfrak{d}} \left[ \Lambda_{\mathfrak{d}}^*(t) - \Lambda_{\mathfrak{d}}^*\left( \frac{1}{2} \right) \right] + \mathfrak{d} \left[ t(1-t) - \frac{1}{4} \right] \right)$$

$$\le \sqrt{\mathfrak{d}} \sup_{|t-1/2| \le c/\mathfrak{d}^{1/4}} \left| \Lambda_{\mathfrak{d}}^*(t) - \Lambda_{\mathfrak{d}}^*\left( \frac{1}{2} \right) \right| + \mathfrak{d} \sup_{0 < t < 1/2} \left[ t(1-t) - \frac{1}{4} \right]$$

$$= o_P(\sqrt{\mathfrak{d}}),$$

which follows from the weak convergence in (B.1) and the continuity of the limit process applied to the first term on the right-hand side of the latter display and obvious reasoning in case of the second term. Relation (B.3) is established. To obtain the assertion of Remark 2.1 in the case of $\Lambda(\mathfrak{d})$, it suffices therefore to prove asymptotic normality for $\sum_{\ell=1}^{\mathfrak{d}} (B_\ell^2(1/2) - 1/4)/(\mathfrak{d}/8)^{1/2}$. This is, however, immediately implied by the CLT.

To prove the assertion of Remark 2.1 for $\Omega(\mathfrak{d})$, we use the relations $E[\int_0^1 B_\ell^2(t)\,dt] = 1/6$ and $\mathrm{Var}(\int_0^1 B_\ell^2(t)\,dt) = 1/45$, and the central limit theorem. The proof is now complete.  □

PROOF OF THEOREM 2.2.  It is an immediate consequence of the ergodic theorem that, as $n \to \infty$,

$$\sup_{t \in [0,\theta]} \left| \frac{1}{n} \sum_{j=1}^{\lfloor nt \rfloor} \mathbf{y}_j \mathbf{y}_j^T - tE[\mathbf{y}_0 \mathbf{y}_0^T] \right| = o_P(1),$$

$$\sup_{t \in (\theta,1]} \left| \frac{1}{n} \sum_{j=\lfloor n\theta \rfloor + 1}^{\lfloor nt \rfloor} \mathbf{y}_j^*(\mathbf{y}_j^*)^T - (t-\theta)E[\mathbf{y}_0^*(\mathbf{y}_0^*)^T] \right| = o_P(1).$$

By assumption (2.11), it follows that

$$\sup_{t \in (\theta,1]} \left| \frac{1}{n} \sum_{j=\lfloor n\theta \rfloor + 1}^{\lfloor nt \rfloor} \mathbf{z}_{j,n} \mathbf{z}_{j,n}^T \right| \leq \frac{1}{n} \sum_{j=\lfloor n\theta \rfloor + 1}^{\lfloor nt \rfloor} |\mathbf{z}_{j,n} \mathbf{z}_{j,n}^T| = o_P(1).$$

Since similar arguments also imply that

$$\sup_{t \in (\theta,1]} \left| \frac{1}{n} \sum_{j=\lfloor n\theta \rfloor + 1}^{\lfloor nt \rfloor} \mathbf{y}_j^* \mathbf{z}_{j,n}^T \right| = \sup_{t \in (\theta,1]} \left| \frac{1}{n} \sum_{j=\lfloor n\theta \rfloor + 1}^{\lfloor nt \rfloor} \mathbf{z}_{j,n}(\mathbf{y}_j^*)^T \right| = o_P(1),$$

the proof of Theorem 2.2 is complete.  □

**B.2. Properties of the change-point estimator.**  The purpose of this subsection is to outline a proof that provides us with the relation in (2.13). In order to keep the arguments short, we simplify the setting and introduce the following technical assumptions. To begin with, we require the pre-change sequence $(\mathbf{y}_j : j \leq k^*)$ to be as in Assumption 2.2, while the post-change sequence $(\mathbf{y}_j : j > k^*)$ is given by $\mathbf{y}_j = \mathbf{y}_j^*$ with $(\mathbf{y}_j^* : j \in \mathbb{Z})$ as in Assumption 2.2, thereby omitting the additional noise terms $\mathbf{z}_{j,n}$ here.

Let $\mathbf{x}_j = \mathrm{vech}[\mathbf{y}_j \mathbf{y}_j^T]$ and $\mathbf{x}_j^* = \mathrm{vech}[\mathbf{y}_j^*(\mathbf{y}_j^*)^T]$ and assume that the thus defined sequences satisfy the strong law of large numbers with a rate $a_n$, that is,

$$\frac{1}{a_n} \sum_{j=1}^n (\mathbf{x}_j - \boldsymbol{\mu}) \to 0 \quad \text{and} \quad \frac{1}{a_n} \sum_{j=1}^n (\mathbf{x}_j^* - \boldsymbol{\mu}^*) \to 0 \qquad \text{a.s.,}$$

where $\boldsymbol{\mu} = E[\mathbf{x}_j]$ and $\boldsymbol{\mu}^* = E[\mathbf{x}_j^*]$. Since we are under $H_A$, we have that $\boldsymbol{\Delta} = \boldsymbol{\mu} - \boldsymbol{\mu}^* \neq 0$.

THEOREM B.1. *Let $k^* = [n\theta]$ for some $\theta \in (0, 1)$. In the above described setting, $\hat{\theta}_n$ is strongly consistent for $\theta$ if*

$$\left(\frac{a_n}{n}\right)^2 \frac{1}{\boldsymbol{\Delta}^T \hat{\Sigma}_n^{-1} \boldsymbol{\Delta}} \to 0 \qquad (n \to \infty),$$

*provided that $\hat{\Sigma}_n \to \Sigma$ with probability one and $\Sigma$ is positive definite.*

If the quantity $\boldsymbol{\Delta}_n$ does not depend on the sample size $n$, then $\boldsymbol{\Delta}^T \hat{\Sigma}_n^{-1} \boldsymbol{\Delta} > 0$ converges to a positive constant and the sufficient condition stated in Theorem B.1 reduces to $a_n/n \to 0$. The latter is essentially satisfied if strong laws of large numbers hold for the pre-change and post-change sequences. To establish Theorem B.1, we need the following lemma which is given without proof.

LEMMA B.1. *Let $f(t)$ and $g(t)$ be functions on $[0, \theta]$ of which $f(t)$ is increasing. As long as $f(\theta) - f(\theta - \gamma) \geq \sup_t |g(t)|$ for some $\gamma \in [0, \theta]$, we have that*

$$\arg\max_{t \in [0,\theta]} [f(t) + g(t)] \geq \theta - \gamma.$$

*An analogous result can be stated if $f(t)$ and $g(t)$ are functions on $[\theta, 1]$ of which $f(t)$ is decreasing.*

PROOF OF THEOREM B.1. Observe first that $\mathcal{L}_n = \boldsymbol{\Delta}^T \hat{\Sigma}_n^{-1} \boldsymbol{\Delta} > 0$ because $\boldsymbol{\Delta} \neq 0$ and $\hat{\Sigma}_n^{-1}$ is positive definite by assumption. Next, straightforward computations reveal that $\mathcal{S}_{nt} \hat{\Sigma}_n^{-1} \mathcal{S}_{nt} = Q_1(nt) + Q_2(nt) + Q_3(nt)$, where

$$Q_1(nt) = \begin{cases} nt^2(1-\theta)^2 \mathcal{L}_n, & t \in [0, \theta], \\ n\theta^2(1-t)^2 \mathcal{L}_n, & t \in (\theta, 1], \end{cases}$$

$$Q_2(nt) = \frac{1}{n} \boldsymbol{\varepsilon}(nt)^T \hat{\Sigma}_n^{-1} \boldsymbol{\varepsilon}(nt), \qquad \boldsymbol{\varepsilon}(nt) = \sum_{j=1}^{[nt]} \boldsymbol{\varepsilon}_j - \frac{[nt]}{n} \sum_{j=1}^{n} \boldsymbol{\varepsilon}_j,$$

and

$$Q_3(nt) = \begin{cases} \dfrac{2t(1-\theta)}{\sqrt{n}} \boldsymbol{\Delta}^T \hat{\Sigma}_n^{-1} \boldsymbol{\varepsilon}(nt), & t \in [0, \theta], \\ \dfrac{2\theta(1-t)}{\sqrt{n}} \boldsymbol{\Delta}^T \hat{\Sigma}_n^{-1} \boldsymbol{\varepsilon}(nt), & t \in (\theta, 1]. \end{cases}$$

The quadratic form $Q_1(nt)$ is clearly increasing in $[0, \theta]$ and decreasing on $[\theta, 1]$ taking its maximum in $\theta$. It remains therefore to show that $Q_2(nt) + Q_3(nt)$ is small compared to $Q_1(nt)$ in the sense of the condition given in Lemma B.1. To this end note that, for vectors $\mathbf{a}, \mathbf{b} \in \mathbb{R}^d$ and a quadratic matrix $M \in \mathbb{R}^{d \times d}$, $|\mathbf{a} M \mathbf{b}| \leq |\mathbf{a}| |M|_E |\mathbf{b}|$. Thus, for all $n \geq n_0$,

$$\sup_{t \in [0,1]} Q_2(nt) \leq \frac{1}{n} \sup_{t \in [0,1]} |\boldsymbol{\varepsilon}(nt)|^2 |\hat{\Sigma}_n^{-1}|_E$$

$$\leq \frac{2}{n} \sup_{t \in [0,1]} |\boldsymbol{\varepsilon}(nt)|^2 |\Sigma^{-1}|_E = o\left(\frac{a_n^2}{n}\right) \qquad \text{a.s.,}$$

since $|\Sigma^{-1}|_E < \infty$ by the positive definiteness of $\Sigma$. Further, for all $n \geq n_0$,

$$\sup_{t \in [0,1]} |Q_3(nt)| \leq \frac{4\theta(1-\theta)}{\sqrt{n}} |\Delta| |\Sigma^{-1}|_E \sup_{t \in [0,1]} |\boldsymbol{\varepsilon}(nt)| = o\left(\frac{\Delta a_n}{\sqrt{n}}\right) \qquad \text{a.s.}$$

Since $a_n/\sqrt{n} \to \infty$ (our sequence satisfies a CLT, therefore the law of large numbers cannot hold with a rate smaller than $\sqrt{n}$), we have shown that $\sup_{t \in [0,1]} |Q_2(nt) + Q_3(nt)| = o(a_n^2/n)$ with probability one. Utilizing Lemma B.1 with $f_n(t) = Q_1(nt)$ and $g_n(t) = Q_2(nt) + Q_3(nt)$ yields that

$$\arg\max_{t \in [0,\theta]} \mathcal{S}_{nt} \hat{\Sigma}_n^{-1} \mathcal{S}_{nt} \geq \theta - \gamma_n,$$

provided $f_n(\theta) - f_n(\theta - \gamma_n) = n\gamma_n(2\theta - \gamma_n)(1-\theta)^2 \mathcal{L}_n \geq \sup_{t \in [0,\theta]} |g_n(t)|$. A similar argument applies also to the inverval $[\theta, 1]$. Choosing $\gamma_n = (a_n/n)^2/\mathcal{L}_n$ verifies the theorem. $\square$

In the special but widely applicable case of weakly dependent sequences satisfying the law of the iterated logarithm, one can pick $a_n = \sqrt{n \log\log n}$. If, in addition, $\boldsymbol{\Delta}$ is also independent of $n$, then the rate given in (2.13) follows.

### B.3. Proofs of Theorems 4.1–4.5.

PROOF OF THEOREM 4.1. We need to verify the assumptions of Theorem A.2. Note first that the moment assumptions are satisfied. Moreover, we can define the approximating random vectors $(\mathbf{y}_j^{(m)} : j \in \mathbb{Z})$ by letting $\mathbf{y}_j^{(m)} = \sum_{\ell=0}^m C_\ell \boldsymbol{\varepsilon}_{j-\ell}$ for $j \in \mathbb{Z}$. Then, condition (2.5) is satisfied, since

$$\sum_{m \geq 1} \|\mathbf{y}_0 - \mathbf{y}_0^{(m)}\|_4 \leq \|\boldsymbol{\varepsilon}_0\|_4 \sum_{m \geq 1} \sum_{\ell \geq m+1} |C_\ell| < \infty$$

by an application of Minkowski's inequality and assumption on the matrices $(C_\ell : \ell \geq 0)$. This is the assertion. $\square$

PROOF OF THEOREM 4.2. Using the Hadamard product $\circ$ in the defining equations (4.3) and (4.4) allows for proving the assertion componentwise.

Let therefore $(y_j : j \in \mathbb{Z})$ denote a generic coordinate of $(\mathbf{y}_j : j \in \mathbb{Z})$. As in the statement of Theorem 4.2, we add—if necessary—zero coefficients and assume that $(y_j : j \in \mathbb{Z})$ follows the GARCH$(r, r)$ model, $r = \max\{p, q\}$,

$$y_j = \sigma_j \varepsilon_j,$$

$$\sigma_j^2 = \omega + \sum_{\ell=1}^{r} \alpha_\ell \sigma_{j-\ell}^2 + \sum_{\ell=1}^{r} \beta_\ell y_{j-\ell}^2.$$

The assumption $\gamma_C < 1$ ensures that $(y_j : j \in \mathbb{Z})$ has a unique strictly stationary solution (see [13, 14]). Berkes, Hörmann and Horváth [8] have shown that this solution has the representation

$$(B.4) \quad y_j = \sqrt{\omega} \varepsilon_j \left( 1 + \sum_{n=1}^{\infty} \sum_{1 \le \ell_1, \ldots, \ell_n \le r} \prod_{i=1}^{n} (\alpha_{\ell_i} + \beta_{\ell_i} \varepsilon_{j-\ell_1-\cdots-\ell_i}^2) \right)^{1/2}.$$

From this representation and Minkowski's inequality, we obtain that

$$(E[y_j^4])^{1/2} \le \omega \|\varepsilon_0\|_4^2 \left( 1 + \sum_{n=1}^{\infty} \sum_{1 \le \ell_1, \ldots, \ell_n \le r} \prod_{i=1}^{n} \|\alpha_{\ell_i} + \beta_{\ell_i} \varepsilon_0^2\|_2 \right)$$

$$\le \omega \|\varepsilon_0\|_4^2 \sum_{n=0}^{\infty} \gamma_C^n$$

and consequently $E[y_j^4] < \infty$ on behalf of $\gamma_C < 1$. It remains to check condition (2.5) coordinatewise. In view of (B.4), a natural candidate for the approximating variables $(y_j^{(m)} : j \in \mathbb{Z})$ is given by the truncated variables

$$(B.5) \quad y_j^{(m)} = \sqrt{\omega} \varepsilon_j \left( 1 + \sum_{n=1}^{[m/r]} \sum_{1 \le \ell_1, \ldots, \ell_n \le r} \prod_{i=1}^{n} (\alpha_{\ell_i} + \beta_{\ell_i} \varepsilon_{j-\ell_1-\cdots-\ell_i}^2) \right)^{1/2}.$$

Observe that $y_j^{(m)}$ is measurable with respect to the $\sigma$-algebra generated by the random variables $\varepsilon_j, \varepsilon_{j-1}, \ldots, \varepsilon_{j-m}$. Using the fact that $|y_j - y_j^{(m)}|^4 \le |y_j^2 - (y_j^{(m)})^2|^2$, it is now easy to see that (2.5) holds if $\gamma_C < 1$. The assertion of Theorem 4.2 is consequently implied by Theorem A.2. $\square$

PROOF OF THEOREM 4.3. At first, we identify the only candidate for the strictly stationary solution iterating the defining volatility equation (4.5). Note that, by (4.3), $\mathbf{y}_j \circ \mathbf{y}_j = \boldsymbol{\sigma}_j \circ \boldsymbol{\sigma}_j \circ \boldsymbol{\varepsilon}_j \circ \boldsymbol{\varepsilon}_j$, and therefore

$$\boldsymbol{\sigma}_j \circ \boldsymbol{\sigma}_j = \boldsymbol{\omega} + \sum_{\ell=1}^{p} A_\ell (\boldsymbol{\sigma}_{j-\ell} \circ \boldsymbol{\sigma}_{j-\ell}) + \sum_{\ell=1}^{q} B_\ell (\mathbf{y}_{j-\ell} \circ \mathbf{y}_{j-\ell})$$

$$(B.6) \quad = \boldsymbol{\omega} + \sum_{\ell=1}^{r} (A_\ell + B_\ell E_{j-\ell}) (\boldsymbol{\sigma}_{j-\ell} \circ \boldsymbol{\sigma}_{j-\ell})$$

$$= \boldsymbol{\omega} + \left[\sum_{k=1}^{\infty} \sum_{1 \le \ell_1, \ldots, \ell_k \le r} \prod_{i=1}^{k}(A_{\ell_i} + B_{\ell_i} E_{j-\ell_1-\cdots-\ell_i})\right] \boldsymbol{\omega},$$

where $r = \max\{p, q\}$, $E_j = \mathrm{diag}(\boldsymbol{\varepsilon}_j \circ \boldsymbol{\varepsilon}_j)$, $j \in \mathbb{Z}$. Utilizing the contraction condition $\gamma_{J,\alpha} < 1$, we will show that the expression on the right-hand side of (B.6) has a finite $\|\cdot\|_{E,\alpha}$-norm. To this end, write

$$\left\|\sum_{k=1}^{\infty} \sum_{1 \le \ell_1, \ldots, \ell_k \le r} \prod_{i=1}^{k}(A_{\ell_i} + B_{\ell_i} E_{j-\ell_1-\cdots-\ell_i})\right\|_{E,\alpha}$$

$$\le \sum_{k=1}^{\infty} \sum_{1 \le \ell_1, \ldots, \ell_k \le r} \prod_{i=1}^{k} \|A_{\ell_i} + B_{\ell_i} E_0\|_{E,\alpha}$$

$$= \sum_{k=1}^{\infty} \left(\sum_{\ell=1}^{r} \|A_\ell + B_\ell E_0\|_{E,\alpha}\right)^k$$

$$= \sum_{k=1}^{\infty} \gamma_{J,\alpha}^k < \infty,$$

where the first inequality sign is obtained after an application of the matrix norm inequality $|MN|_E \le |M|_E |N|_E$, Minkowski's inequality for the $L^\alpha$-norm and the fact that $(E_j : j \in \mathbb{Z})$ are independent, identically distributed random matrices. The proof is complete. $\square$

PROOF OF THEOREM 4.4. We need to verify Assumption 2.1. Condition (2.3) has been established in Theorem 4.3. It remains to determine an approximating sequence $(\mathbf{y}_j^{(m)} : j \in \mathbb{Z})$ that satisfies (2.5). To this end, introduce

$$\mathbf{y}_0^{(m)} = \boldsymbol{\sigma}_0^{(m)} \circ \boldsymbol{\varepsilon}_0,$$

$$\boldsymbol{\sigma}_0^{(m)} \circ \boldsymbol{\sigma}_0^{(m)} = \boldsymbol{\omega} + \left[\sum_{k=1}^{m} \sum_{1 \le \ell_1, \ldots, \ell_k \le r} \prod_{i=1}^{k}(A_{\ell_i} + B_{\ell_i} E_{-\ell_1-\cdots-\ell_i})\right] \boldsymbol{\omega},$$

where $r = \max\{p, q\}$. Due to the stationarity result given in Theorem 4.3, it suffices to prove (2.5) for $j = 0$. This will be achieved in three steps. We first establish a result for the distance between $\boldsymbol{\sigma}_0$ and the approximating volatility $\boldsymbol{\sigma}_0^{(m)}$. Observe that, using the definition of the Euclidean norm and the $L^\alpha$-norm, we obtain that

$$\|\boldsymbol{\sigma}_0 - \boldsymbol{\sigma}_0^{(m)}\|_4^2 = \left\|\sum_{i=1}^{d}(\sigma_{0,i} - \sigma_{0,i}^{(m)})^2\right\|_2$$

$$\leq \sum_{i=1}^{d} \|(\sigma_{0,i} - \sigma_{0,i}^{(m)})^2\|_2$$

$$\leq \sum_{i=1}^{d} \|\sigma_{0,i}^2 - (\sigma_{0,i}^{(m)})^2\|_2$$

$$(B.7) \qquad \leq d \max_{1 \leq i \leq d} \|\sigma_{0,i}^2 - (\sigma_{0,i}^{(m)})^2\|_2$$

$$= d\Big(\max_{1 \leq i \leq d} \|\sigma_{0,i}^2 - (\sigma_{0,i}^{(m)})^2\|_2^2\Big)^{1/2}$$

$$\leq d\Big(\sum_{i=1}^{d} \|\sigma_{0,i}^2 - (\sigma_{0,i}^{(m)})^2\|_2^2\Big)^{1/2}$$

$$= d\|\boldsymbol{\sigma}_0 \circ \boldsymbol{\sigma}_0 - \boldsymbol{\sigma}_0^{(m)} \circ \boldsymbol{\sigma}_0^{(m)}\|_2,$$

where we have applied that $|a - b|^2 \leq |a^2 - b^2|$ if $a, b \geq 0$ to obtain the first inequality sign in the second line of the display. We further estimate the distance between the Hadamard products of the volatilities $\boldsymbol{\sigma}_0 \circ \boldsymbol{\sigma}_0$ and $\boldsymbol{\sigma}_0^{(m)} \circ \boldsymbol{\sigma}_0^{(m)}$. We have

$$\|\boldsymbol{\sigma}_0 \circ \boldsymbol{\sigma}_0 - \boldsymbol{\sigma}_0^{(m)} \circ \boldsymbol{\sigma}_0^{(m)}\|_2$$

$$\leq |\boldsymbol{\omega}| \Big\| \sum_{k=m+1}^{\infty} \sum_{1 \leq \ell_1, \ldots, \ell_k \leq r} \prod_{i=1}^{k} (A_{\ell_i} + B_{\ell_i} E_{-\ell_1 - \cdots - \ell_i}) \Big\|_2$$

$$(B.8) \qquad \leq |\boldsymbol{\omega}| \sum_{k=m+1}^{\infty} \sum_{1 \leq \ell_1, \ldots, \ell_k \leq r} \prod_{i=1}^{k} \|A_{\ell_i} + B_{\ell_i} E_0\|_2$$

$$= |\boldsymbol{\omega}| \sum_{k=m+1}^{\infty} \Big(\sum_{\ell=1}^{r} \|A_\ell + B_\ell E_0\|_2\Big)^k = \frac{|\boldsymbol{\omega}| \gamma_J^{m+1}}{1 - \gamma_J}.$$

Finally, combining (B.7) with (B.8), we arrive at

$$\|\mathbf{y}_0 - \mathbf{y}_0^{(m)}\|_4 \leq \|\boldsymbol{\varepsilon}_0\|_4 \|\boldsymbol{\sigma}_0 - \boldsymbol{\sigma}_0^{(m)}\|_4$$

$$\leq \sqrt{d} \|\boldsymbol{\varepsilon}_0\|_4 (\|\boldsymbol{\sigma}_0 \circ \boldsymbol{\sigma}_0 - \boldsymbol{\sigma}_0^{(m)} \circ \boldsymbol{\sigma}_0^{(m)}\|_2)^{1/2}$$

$$\leq \sqrt{d} \|\boldsymbol{\varepsilon}_0\|_4 \Big(\frac{|\boldsymbol{\omega}| \gamma_J^{m+1}}{1 - \gamma_J}\Big)^{1/2}.$$

The right-hand side is summable on account of $\gamma_J < 1$, thus yielding (2.5) and completing the proof of Theorem 4.4. $\square$

PROOF OF THEOREM 4.5. It follows immediately from Theorems 4.3 and 4.4. $\square$

**B.4. Proofs of Theorems 4.6 and 4.7.**

PROOF OF THEOREM 4.6. The assertion is an immediate consequence of Brandt [15] and Bougerol and Picard [13]. □

The proof of Theorem 4.7 requires three auxiliary lemmas and some basic results from linear algebra. The first lemma establishes connections between the various norms that are in use.

LEMMA B.2. *If $M$ is a $d \times d$ matrix, then it holds that $|M|_E \leq \sqrt{2} |\operatorname{vech}[M]|$.*

PROOF. The inequality follows from the corresponding norm properties. □

LEMMA B.3. *Let $M$ and $N$ be two symmetric square matrices, and let*

$$\exp(M) = \sum_{k=0}^{\infty} \frac{1}{k!} M^k \quad and \quad \exp(N) = \sum_{k=0}^{\infty} \frac{1}{k!} N^k.$$

*Then, it holds that:*

  (i) *$\exp(M)$ is positive definite;*
  (ii) *for every real numbers $a$ and $b$, $\exp(aM)\exp(bM) = \exp((a+b)M)$;*
  (iii) *$\exp(M)\exp(-M) = I$, where $I$ denotes the identity matrix;*
  (iv) *$(\exp(M))^{1/2} = \exp(\frac{1}{2}M)$;*
  (v) *$|\exp(M + N) - \exp(M)|_E \leq |N|_E \exp(|M|_E) \exp(|N|_E)$.*

PROOF. See, for example, Horn and Johnson [30]. □

LEMMA B.4. *Let $X$ be a positive random variable such that $E[\exp(tX)] < \infty$ if $t < \tau$. Then, there is a $t_0 \in (0, \tau)$ such that $E[\exp(tX)] \leq 1 + 2tE[X]$ for all $t \in [0, t_0]$.*

PROOF. The assertion follows from routine arguments. □

PROOF OF THEOREM 4.7. To verify Assumption 2.1 note that (2.3) follows from Theorem 4.6. It remains to define a suitable sequence $(\mathbf{y}_j^{(m)} : j \in \mathbb{Z})$ that satisfies (2.5). To this end, let

$$\mathbf{y}_0^{(m)} = (H_0^{(m)})^{1/2} \boldsymbol{\varepsilon}_0,$$

where the matrix $H_0^{(m)}$ is defined via truncating the strictly stationary solution at lag $m$ as

$$\operatorname{vech}[\log H_0^{(m)} - C] = \sum_{k=0}^{m} A^k F(\boldsymbol{\varepsilon}_{-k-1}, \ldots, \boldsymbol{\varepsilon}_{-k-q}).$$

Next, let

$$S_0 = C + \text{math}\left[\sum_{k=0}^{\infty} A^k F(\boldsymbol{\varepsilon}_{-k-1}, \ldots, \boldsymbol{\varepsilon}_{-k-q})\right],$$

$$S_0^{(m)} = C + \text{math}\left[\sum_{k=0}^{m} A^k F(\boldsymbol{\varepsilon}_{-k-1}, \ldots, \boldsymbol{\varepsilon}_{-k-q})\right],$$

where math denotes the inverse operator of vech. Then, by assumption on $C$, the matrices $S_0$ and $S_0^{(m)}$ are symmetric. It holds, moreover, that $H_0 = \exp(S_0)$, $H_0^{(m)} = \exp(S_0^{(m)})$, $H_0^{1/2} = \exp(\frac{1}{2}S_0)$ and $(H_0^{(m)})^{1/2} = \exp(\frac{1}{2}S_0^{(m)})$. For the last two statements, part (iv) of Lemma B.3 has been applied. Utilizing the previous statements, we estimate next the matrix norm difference between $H_0^{1/2}$ and $(H_0^{(m)})^{1/2}$ as

$$|H_0^{1/2} - (H_0^{(m)})^{1/2}|_E = |\exp(\tfrac{1}{2}[S_0 + S_0^{(m)} - S_0^{(m)}]) - \exp(\tfrac{1}{2}S_0^{(m)})|_E$$
$$\leq |S_0 - S_0^{(m)}|_E \exp(\tfrac{1}{2}|S_0^{(m)}|_E) \exp(\tfrac{1}{2}|S_0 - S_0^{(m)}|_E),$$

where we have applied part (v) of Lemma B.3 with $M = \frac{1}{2}S_0^{(m)}$ and $N = \frac{1}{2}(S_0 - S_0^{(m)})$ to obtain the inequality sign (additionally dropping a factor $\frac{1}{2}$ on the right-hand side). Since $\boldsymbol{\varepsilon}_0$ is independent of $S_0$ and $S_0^{(m)}$, we conclude further that

$$\begin{aligned}
&\|\mathbf{y}_0 - \mathbf{y}_0^{(m)}\|_4 \\
&\quad \leq \|\boldsymbol{\varepsilon}_0\|_4 \|H_0^{1/2} - (H_0^{(m)})^{1/2}\|_{E,4} \\
\text{(B.9)} \quad &\quad \leq \|\boldsymbol{\varepsilon}_0\|_4 \||S_0 - S_0^{(m)}|_E \exp(\tfrac{1}{2}|S_0^{(m)}|_E) \exp(\tfrac{1}{2}|S_0 - S_0^{(m)}|_E)\|_4 \\
&\quad \leq \|\boldsymbol{\varepsilon}_0\|_4 \|S_0 - S_0^{(m)}\|_{E,8\beta} \|\exp(\tfrac{1}{2}|S_0^{(m)}|_E)\|_{4\alpha} \\
&\qquad \times \|\exp(\tfrac{1}{2}|S_0 - S_0^{(m)}|_E)\|_{8\beta}.
\end{aligned}$$

The latter inequality sign follows from an application of Hölder's inequality with $\alpha, \beta > 1$ such that $1/\alpha + 1/\beta = 1$. In the following, we investigate the norms in (B.9). By assumption, $\|\boldsymbol{\varepsilon}_0\|_4 < \infty$. Lemma B.2, Minkowski's inequality and the existence of a moment generating function for $F(\boldsymbol{\varepsilon}_{-1}, \ldots, \boldsymbol{\varepsilon}_{-q})$ imply that

$$\|S_0 - S_0^{(m)}\|_{E,8\beta} \leq \sqrt{2} \|\text{vech}[S_0 - S_0^{(m)}]\|_{8\beta}$$
$$\leq \sqrt{2} \left\|\sum_{k=m+1}^{\infty} |A|_E^k |F(\boldsymbol{\varepsilon}_{-k-1}, \ldots, \boldsymbol{\varepsilon}_{-k-q})|\right\|_{8\beta}$$

$$\leq \sqrt{2} \sum_{k=m+1}^{\infty} |A|_E^k \|F(\varepsilon_{-k-1}, \ldots, \varepsilon_{-k-q})\|_{8\beta}$$

$$\leq C|A|_E^m < \infty,$$

where $C > 0$ is a constant, on account of $|A|_E < 1$ by assumption. By similar arguments, we continue estimating

$$\left\| \exp\left( \frac{1}{2} |S_0 - S_0^{(m)}|_E \right) \right\|_{8\beta}$$

$$\leq \left\| \exp\left( \frac{1}{\sqrt{2}} |\mathrm{vech}[S_0 - S_0^{(m)}]| \right) \right\|_{8\beta}$$

$$\leq \left\| \exp\left( \frac{1}{\sqrt{2}} \left| \sum_{k=m+1}^{\infty} A^k F(\varepsilon_{-k-1}, \ldots, \varepsilon_{-k-q}) \right| \right) \right\|_{8\beta}$$

$$\leq \left\| \exp\left( \frac{1}{\sqrt{2}} \sum_{k=m+1}^{\infty} |A|_E^k |F(\varepsilon_{-k-1}, \ldots, \varepsilon_{-k-q})| \right) \right\|_{8\beta}$$

$$= \mathcal{I}(m, \beta).$$

Observe that $(F(\varepsilon_{-k-1}, \ldots, \varepsilon_{-k-q}) : k \geq 0)$ is a sequence of $q$-dependent random variables. Therefore, each of the sequences $(F(\varepsilon_{-k-1}, \ldots, \varepsilon_{-k-q}) : k = \nu \bmod q)$, $\nu = 0, 1, \ldots, q-1$, consists of independent random variables. Defining

$$I_\nu(m) = \sum_{m < k = \nu \bmod q} |A|_E^k |F(\varepsilon_{-k-1}, \ldots, \varepsilon_{-k-q})|, \qquad \nu = 0, 1, \ldots, q-1,$$

we obtain that

$$\mathcal{I}(m, \beta) = \left( E\left[ \prod_{\nu=0}^{q-1} \exp\left( \frac{8\beta}{\sqrt{2}} I_\nu(m) \right) \right] \right)^{1/(8\beta)}$$

$$\leq \left( \prod_{\nu=0}^{q-1} E\left[ \exp\left( \frac{8\beta q}{\sqrt{2}} I_\nu(m) \right) \right] \right)^{1/(8\beta q)}.$$

It remains to estimate the expectations in the latter equation. We do so first assuming that $m$ is sufficiently large. In this case, the fact that $I_\nu(m)$ is a sum of independent random variables and Lemma B.4 imply that there is a constant $c$ such that

$$E\left[ \exp\left( \frac{8\beta q}{\sqrt{2}} I_\nu(m) \right) \right]$$

$$= \prod_{m < k = \nu \bmod q} E\left[ \exp\left( \frac{8\beta q}{\sqrt{2}} |A|_E^k |F(\varepsilon_{-1}, \ldots, \varepsilon_{-q})| \right) \right]$$

$$\leq \prod_{m < k = \nu \bmod q} (1 + c|A|_E^k E[|F(\boldsymbol{\varepsilon}_{-1}, \ldots, \boldsymbol{\varepsilon}_{-q})|])$$

$$\leq \prod_{m < k = \nu \bmod q} \exp(c|A|_E^k E[|F(\boldsymbol{\varepsilon}_{-1}, \ldots, \boldsymbol{\varepsilon}_{-q})|])$$

$$< \exp(\tilde{c}|A|_E^m),$$

if $m \geq m_0$ for a suitable constant $\tilde{c} > 0$. If $m < m_0$, we set $H_0^{(m)} = 0$. Then the proof is complete if we show that

$$\sup_{1 \leq m \leq \infty} \left\| \exp\left( \frac{1}{2} |S_0^{(m)}|_E \right) \right\|_{4\alpha} < \infty,$$

where $S_0^{(\infty)} = S_0$. This, however, follows from similar arguments as before as long as the moment generating function $E[\exp(t|F(\boldsymbol{\varepsilon}_{-1}, \ldots, \boldsymbol{\varepsilon}_{-q})|)]$ is finite for $t \geq 4\alpha q/\sqrt{2}$. Now pick $\alpha$ close to 1 to obtain the assertion. $\square$

## REFERENCES

[1] ANDREOU, E. and GHYSELS, E. (2002). Detecting multiple breaks in financial market volatility dynamics. *J. Appl. Econometrics* **17** 579–600.

[2] ANDREOU, E. and GHYSELS, E. (2003). Tests for breaks in the conditional co-movements of asset returns. *Statist. Sinica* **13** 1045–1073. MR2026061

[3] ANDREWS, D. W. K. (1984). Nonstrong mixing autoregressive processes. *J. Appl. Probab.* **21** 930–934. MR0766830

[4] ANDREWS, D. W. K. (1991). Heteroskedasticity and autocorrelation-consistent covariance matrix estimation. *Econometrica* **59** 817–858. MR1106513

[5] AUE, A., BERKES, I. and HORVÁTH, L. (2006). Strong approximation for the sums of squares of augmented GARCH sequences. *Bernoulli* **12** 583–608. MR2248229

[6] BAI, J. and PERRON, P. (1998). Estimating and testing linear models with multiple structural changes. *Econometrica* **66** 47–78. MR1616121

[7] BAUWENS, L., LAURENT, S. and ROMBOUTS, J. V. K. (2006). Multivariate GARCH models: A survey. *J. Appl. Econometrics* **21** 79–109. MR2225523

[8] BERKES, I., HÖRMANN, S. and HORVÁTH, L. (2008). The functional central limit theorem for a family of GARCH observations with applications. *Statist. Probab. Lett.* **78** 2725–2730. MR2465114

[9] BERKES, I., HORVÁTH, L. and KOKOSZKA, P. (2003). GARCH processes: Structure and estimation. *Bernoulli* **9** 201–227. MR1997027

[10] BILLINGSLEY, P. (1968). *Convergence of Probability Measures.* Wiley, New York. MR0233396

[11] BOLLERSLEV, T. (1986). Generalized autoregressive conditional heteroskedasticity. *J. Econometrics* **31** 307–327. MR0853051

[12] BOLLERSLEV, T. (1990). Modeling the coherence in short-run nominal exchange rates: A multivariate generalized ARCH model. *Rev. Econom. Statist.* **74** 498–505.

[13] BOUGEROL, P. and PICARD, N. (1992). Strict stationarity of generalized autoregressive processes. *Ann. Probab.* **20** 1714–1730. MR1188039

[14] BOUGEROL, P. and PICARD, N. (1992). Stationarity of GARCH processes and of some nonnegative time series. *J. Econometrics* **52** 115–127. MR1165646

[15] BRANDT, A. (1986). The stochastic equation $Y_{n+1} = A_n Y_n + B_n$ with stationary coefficients. *Adv. in Appl. Probab.* **18** 211–220. MR0827336

[16] BROCKWELL, P. J. and DAVIS, R. A. (1991). *Time Series: Theory and Methods*, 2nd ed. Springer, New York. MR1093459

[17] CARRASCO, M. and CHEN, X. (2002). Mixing and moment properties of various GARCH and stochastic volatility models. *Econometric Theory* **18** 17–39. MR1885348

[18] CSÖRGŐ, M. and HORVÁTH, L. (1997). *Limit Theorems in Change-Point Analysis*. Wiley, Chichester.

[19] DAVIDSON, J. (1994). *Stochastic Limit Theory*. Oxford Univ. Press, Oxford. MR1430804

[20] DAVIS, R. A., LEE, T. C. M. and RODRIGUEZ-YAM, G. (2008). Break detection for a class of nonlinear time series models. *J. Time Ser. Anal.* **29** 834–867. MR2450899

[21] DAVIS, R. A. and MIKOSCH, T. (1998). The sample autocorrelations of heavy-tailed processes with applications to ARCH. *Ann. Statist.* **26** 2049–2080. MR1673289

[22] DOUKHAN, P. and LOUHICHI, S. (1999). A new weak dependence condition and application to moment inequalities. *Stochastic Process. Appl.* **84** 313–343. MR1719345

[23] ENGLE, R. F. (1982). Autoregressive conditional heteroskedasticity with estimates of the variance of U.K. inflation. *Econometrica* **50** 987–1008. MR0666121

[24] ENGLE, R. F., NG, V. K. and ROTHSCHILD, M. (1990). Asset pricing with a factor ARCH covariance structure: Empirical estimates for treasury bills. *J. Econometrics* **45** 213–238.

[25] ERDÉLYI, A. (1954). *Tables of Integral Transforms* **1**. McGraw-Hill, New York.

[26] GOMBAY, E., HORVÁTH, L. and HUŠKOVÁ, M. (1996). Estimators and tests for change in the variance. *Statistics and Decisions* **14** 145–159. MR1406521

[27] HALL, P. and YAO, Q. (2003). Inference in ARCH and GARCH models with heavy-tailed errors. *Econometrica* **71** 285–317. MR1956860

[28] HE, C. and TERÄSVIRTA, T. (1999). Fourth moment structure of the Garch($p, q$) process. *Econometric Theory* **15** 824–846.

[29] HÖRMANN, S. (2008). Augmented GARCH sequences: Dependence structure and asymptotics. *Bernoulli* **14** 543–561.

[30] HORN, R. A. and JOHNSON, C. R. (1991). *Topics in Matrix Analysis*. Cambridge Univ. Press, Cambridge. MR1091716

[31] HSU, D. A. (1979). Detecting shifts of parameter in gamma sequences with applications to stock price and air traffic flow analysis. *J. Amer. Statist. Assoc.* **74** 31–40.

[32] IBRAGIMOV, I. A. (1962). Some limit theorems for stationary processes. *Theory Probab. Appl.* **7** 349–382. MR0148125

[33] INCLÁN, C. and TIAO, G. C. (1994). Use of cumulative sums of squares for retrospective detection of changes of variance. *J. Amer. Statist. Assoc.* **89** 913–923. MR1294735

[34] JEANTHEAU, T. (1998). Strong consistency of estimators for multivariate ARCH models. *Econometric Theory* **14** 70–86. MR1613694

[35] KAWAKATSU, H. (2006). Matrix exponential GARCH. *J. Econometrics* **134** 95–128. MR2328317

[36] KIEFER, J. (1959). $K$-sample analogues of the Kolmogorov–Smirnov and Cramér–v. Mises tests. *Ann. Math. Statist.* **30** 420–447. MR0102882

[37] KOKOSZKA, P. and LEIPUS, R. (2000). Change-point estimation in ARCH models. *Bernoulli* **6** 513–539. MR1762558

[38] LANNE, M. and SAIKKONEN, P. (2007). A multivariate generalized orthogonal factor GARCH model. *J. Bus. Econom. Statist.* **25** 61–75. MR2338871

[39] LI, W. K. (2004). *Diagnostic Checks in Time Series.* Chapman & Hall, Boca Raton, FL.

[40] MIKOSCH, T. (2003). Modeling dependence and tails of financial time series. In: *Extreme Values in Finance, Telecommunications, and the Environment* (B. Finkenstaedt and H. Rootzen, eds.) 185–286. Chapman & Hall, Boca Raton, FL.

[41] MIKOSCH, T. and STRAUMANN, D. (2006). Stable limits of martingale transforms with application to the estimation of GARCH parameters. *Ann. Statist.* **34** 493–522. MR2275251

[42] NELSON, D. B. (1991). Conditional heteroskedasticity in asset returns: A new approach. *Econometrica* **59** 347–370. MR1097532

[43] NZE, P. A. and DOUKHAN, P. (2004). Weak dependence: Models and applications to econometrics. *Econometric Theory* **20** 995–1045. MR2101950

[44] PETROV, V. V. (1995). *Limit Theorems of Probability Theory.* Oxford Univ. Press, Oxford. MR1353441

[45] PÖTSCHER, B. M. and PRUCHA, I. R. (1997). *Dynamic Nonlinear Econometric Models, Asymptotic Theory.* Springer, New York. MR1468737

[46] SERBAN, M., BROCKWELL, A., LEHOCZKY, J. and SRIVASTAVA, S. (2007). Modelling the dynamic dependence structure in multivariate financial time series. *J. Time Ser. Anal.* **28** 763–782. MR2395913

[47] SILVENNOINEN, A. and TERÄSVIRTA, T. (2008). Multivariate GARCH models. In: *Handbook of Financial Time Series* (T. G. Andersen, R. A. Davis, J.-P. Kreiss and T. Mikosch, eds.) 201–229. Springer, New York.

[48] TSAY, R. S. (2005). *Analysis of Financial Time Series*, 2nd ed. Wiley, New York. MR2162112

[49] VRONTOS, I. D., DELLAPORTAS, P. and POLITIS, D. N. (2003). A full-factor multivariate GARCH model. *Econom. J.* **6** 312–334. MR2028238

[50] WU, W. B. (2005). Nonlinear system theory: Another look at dependence. *Proc. Natl. Acad. Sci. USA* **102** 14150–14154. MR2172215

A. AUE
DEPARTMENT OF STATISTICS
UNIVERSITY OF CALIFORNIA
ONE SHIELDS AVENUE
DAVIS, CALIFORNIA 95616
USA
E-MAIL: alexaue@wald.ucdavis.edu

S. HÖRMANN
L. HORVÁTH
DEPARTMENT OF MATHEMATICS
UNIVERSITY OF UTAH
155 SOUTH 1440 EAST
SALT LAKE CITY, UTAH 84112-0090
USA
E-MAIL: hormann@math.utah.edu
horvath@math.utah.edu

M. REIMHERR
DEPARTMENT OF STATISTICS
UNIVERSITY OF CHICAGO
5734 STREET UNIVERSITY AVENUE
CHICAGO, ILLINOIS 60637
USA
E-MAIL: mreimherr@uchicago.edu